\documentclass[11pt]{article}

\usepackage{pinlabel}
\usepackage[all]{xy}

\usepackage{latexsym,epsfig}
\usepackage{amssymb,amsmath,amsfonts}
\usepackage{exscale}
\usepackage{ifthen}
\usepackage{longtable}
\usepackage{subfigure}
\usepackage{caption}
\usepackage{cite}
\usepackage{lscape}

\newcommand{\hl}{\sl}
\newcommand{\hd}{\sl}
 \newcommand{\NN}{\mathbb{N}}
 \newcommand{\ZZ}{\mathbb{Z}}
 \newcommand{\QQ}{\mathbb{Q}}
 \newcommand{\RR}{\mathbb{R}}
\newcommand{\eps}{\varepsilon}
 \newcommand{\XX}{X}
\newcommand{\ganz}{\overline{\XX}}
\newcommand{\rand}{\partial\XX}
 
\newcommand{\regrand}{\partial\XX^{reg}}
\newcommand{\singrand}{\partial\XX^{sing}}
\newcommand{\Frand}{\partial^F\XX}
\newcommand{\Lim}{L_\Gamma}          
\newcommand{\Flim}{F_\Gamma}

\newcommand{\horinf}{\mbox{Vis}^{\infty}}
\newcommand{\horF}{\mbox{Vis}^F}
\newcommand{\rel}{\mbox{Rel}_\Gamma}
\newcommand{\preim}{(\pi^F)^{-1}}
\newcommand{\xo}{{o}}
\newcommand{\Ax}{\mbox{Ax}}
\newcommand{\bs}{{\cal B}}
\newcommand{\pr}{\mbox{pr}}
\newcommand{\is}{\mbox{Is}}
\newcommand{\at}{\!\cdot\!}
\newcommand{\st}{\mbox{such}\ \mbox{that}\ }
\newcommand{\an}{\ \mbox{and}\ }
\newcommand{\wrt}{\mbox{with}\ \mbox{respect}\ \mbox{to}\ }
\newcommand{\be}{\begin{eqnarray*}}
\newcommand{\ee}{\end{eqnarray*}}
\newcommand{\ba}{\begin{array}}
\newcommand{\ea}{\end{array}}

\newcommand{\prf}{{\noindent \sl Proof.\  }}
\newcommand{\qed}{$\hfill\Box$}

\newtheorem{df}{\sc Definition}[section]
\newtheorem{thr}[df]{\bf Theorem}
\newtheorem{lem}[df]{\sc Lemma}
\newtheorem{prp}[df]{\sc Proposition}
\newtheorem{cor}[df]{\sc Corollary}

\begin{document}

\title{Asymptotic geometry in products of Hadamard spaces with rank one isometries}
\author{\sc Gabriele Link\thanks{supported by the FNS grant PP002-102765}}
\date{}
\maketitle

\bibliographystyle{gtart}
\protect\pagenumbering{arabic}
\setcounter{page}{1}

\begin{abstract}
In this article we study asymptotic properties of certain discrete groups $\Gamma$ 
acting by isometries  on a product  $\XX=\XX_1\times \XX_2$ of locally compact Hadamard spaces which admit a geodesic without flat half-plane.  The motivation comes from the fact that  Kac-Moody groups over finite fields, which  can be seen as generalizations of arithmetic groups over function fields, belong to the considered class of groups. Hence one may ask whether classical properties of discrete subgroups of higher rank Lie groups as in  \cite{MR1437472} and \cite{MR1933790} hold in this context. 

In the  first part of the paper we describe the structure of the geometric limit set of $\Gamma$ and prove statements analogous to the results of Benoist in  \cite{MR1437472}. 
The second part  is concerned with the exponential growth rate $\delta_\theta(\Gamma)$ of orbit points in $\XX$  with a prescribed so-called "slope" $\theta\in (0,\pi/2)$, which appropriately generalizes the critical exponent in higher rank.  In analogy to Quint's result in  \cite{MR1933790} we show that  the  homogeneous extension $\Psi_\Gamma$ to $\RR_{\ge 0}^2$ of  $\delta_\theta(\Gamma)$ as a function of $\theta$ is upper semi-continuous and concave.  
\end{abstract}

\maketitle

\section{Introduction}

Let $(\XX_1,d_1)$, $(\XX_2,d_2)$ be Hadamard spaces, i.e. complete simply connected metric spaces of non-positive Alexandrov curvature, and $(\XX,d)$ the product $\XX_1\times \XX_2$ endowed with the metric $d=\sqrt{d_1^2+d_2^2}$. Assume moreover that $\XX_1$, $\XX_2$ are locally compact. 
Each metric space $\XX, \XX_1,\XX_2 $ can be compactified by adding its  geometric boundary $\rand$, $
\rand_1$, $\rand_2$   endowed with the cone topology (see \cite[chapter II]{MR1377265}). It is well-known that the {\hl regular geometric boundary} $\regrand$ of $\XX$ -- which consists of the set of equivalence classes of geodesic rays which do not project to a point in one of the factors -- is a dense open subset of $\rand$ homeomorphic to $\rand_1\times \rand_2\times (0,\pi/2)$. The last factor in this product is called the {\hl slope} of a point in $\regrand$.  The singular geometric boundary $\singrand=\rand\setminus\regrand$ consists of two strata homeomorphic to $\rand_1$, $\rand_2$ respectively. We assign slope $0$ to the first and slope $\pi/2$ to the second one. 

For a group $\Gamma\subseteq \is(\XX_1)\times \is(\XX_2)$ acting properly discontinuously by isometries on $\XX$ we study the limit set $\Lim:=\overline{\Gamma\at x}\cap\rand$, where $x\in\XX$ is arbitrary. Unlike in the case of  CAT$(-1)$-spaces, this geometric limit set is not necessarily a minimal set for the action of $\Gamma$ because an element of $\is(\XX_1)\times \is(\XX_2)$ cannot change the slope $\theta$ of a point in $\rand$. This is similar to the situation in  symmetric spaces or Bruhat-Tits buildings of higher rank. So  by abuse of notation we are going to call the set  $\Frand:=\rand_1\times\rand_2$  the {\hl Furstenberg boundary}, and the projection of $\Lim\cap\regrand$ to $\rand_1\times\rand_2$  the {\hl Furstenberg limit set} $\Flim$ of $\Gamma$. 

In this note we restrict our attention to discrete groups  $\Gamma\subseteq\is(\XX_1)\times\is(\XX_2)$ which contain an element projecting to a rank one element in each factor, i.e. $\Gamma$ contains an element $h=(h_1,h_2)$ \st the invariant geodesics of $h_1,  h_2$ do not bound a flat half-plane in $\XX_1, \XX_2$. Such an isometry of  $\XX$  will be called {\hl regular axial} in the sequel. Moreover, for Theorems A and B below we require as in \cite{DalboKim} that for $i=1,2$ the projection $\Gamma_i$ of $\Gamma$ to $\is(\XX_i)$ is strongly non-elementary: This means that $\Gamma_i$, $i=1,2$, possesses infinitely many limit points and does not globally fix a point at infinity. By Proposition~3.4 in \cite{GR/0809.0470} this condition is equivalent to the fact that both $\Gamma_1$ and $\Gamma_2$ contain a pair of independent rank one elements. For Theorems~C, D and E below we need a slightly stronger assumption: We require that $\Gamma\subseteq\is(\XX_1)\times\is(\XX_2)$ contains  two regular axial isometries  $g=(g_1,g_2)$ and $h=(h_1,h_2)$ \st $g_1, h_1$ and $g_2,h_2$ are pairs of independent rank one elements in $\is(\XX_1)$ resp. $\is(\XX_2)$. 

One important class of examples satisfying our stronger assumption are Kac-Moody groups  $\Gamma$ over a finite field which act by isometries on a product $\XX=\XX_1\times\XX_2$, the CAT$(0)$-realization of the associated twin building ${\cal B}_+\times{\cal B}_-$.  Indeed, there  exists an element $h=(h_1,h_2)$ projecting to a rank one element in each factor by  Remark 5.4 and the proof of Corollary~1.3 in \cite{GR/0809.0470}. Moreover, the action of the Weyl group produces many regular axial isometries $g=(g_1,g_2)$ with $g_i$ independent from $h_i$ for $i=1,2$. Notice that if the order of the ground field is sufficiently large, then $\Gamma\subseteq\is(\XX_1)\times\is(\XX_2)$ is an irreducible lattice (see e.g. \cite{MR1715140} and \cite{GR/0607664}). 

A second type of examples are groups acting properly discontinuously on a product of locally compact 
Hadamard spaces of strictly negative Alexandrov curvature (compare \cite{DalboKim} in the manifold setting). In this special case every non-elliptic and non-parabolic isometry in one of the factors is already a rank one element.  Prominent examples here which are already covered by the results of Y.~Benoist and J.-F.~Quint are Hilbert modular groups acting as irreducible lattices on a product of hyperbolic planes, and graphs of convex cocompact groups of rank one symmetric spaces (see also \cite{MR1230298}). 
But our context is much more general and possible factors include  locally finite, not necessarily regular trees and  Riemannian universal covers of geometric rank one manifolds. 

Our first result is\\[2mm]
{\bf Theorem A}$\quad$ {\sl The Furstenberg limit set is minimal, i.e. $\Flim$ is the smallest non-empty, $\Gamma$-invariant closed subset of $\Frand$.
}\\[3mm]
Moreover we have -- as in the case of symmetric spaces or Bruhat-Tits buildings of higher rank -- the following structure theorem.\\[2mm]
{\bf Theorem B}$\quad$ {\sl The regular geometric limit set splits as a product $\Flim\times P_\Gamma$, where $P_\Gamma\subseteq (0,\pi/2)$ denotes the set of  slopes of  regular limit points. }\\[3mm]
From here on we will assume that $\Gamma$ contains two  regular axial isometries projecting to independent rank one elements in each factor.  Let  $g_1\in \is(\XX_1)$, $g_2\in\is(\XX_2)$  be rank one elements. For $i=1,2$ we denote $g_i^+$ the attractive,  $g_i^-$ the repulsive fixed point, and $l_i(g_i)$ the translation length, i.e. the minimum of the set $\{d_i(x_i,g_ix_i): x_i\in\XX_i\}$. If $g=(g_1,g_2)$, we put $g^+:=(g_1^+,g_2^+)$, $g^-:=(g_1^-,g_2^-)\in\Frand$. Then we have the following two statements:\\[2mm] 
{\bf Theorem C}$\quad$ {\sl $P_\Gamma$ is an interval and we have
$$P_\Gamma=\overline{\{\arctan\big(l_2(g_2)/l_1(g_1)\big): (g_1,g_2)\in\Gamma\,,\ g_1,g_2\ \mbox{rank one}\}} \cap (0,\pi/2)\,.$$ }\\[-4mm]
{\bf Theorem D}$\quad$ {\sl The set of pairs of fixed points $(g^+,g^-)\subset\Frand\times \Frand$ of regular axial isometries in $\Gamma$ is dense in $\big(F_\Gamma
\times F_\Gamma\big)\setminus \Delta$, where $\Delta$ denotes the set of points $(\xi,\eta)$ \st $\xi_1=\eta_1$ or $\xi_2=\eta_2$.   }\\[-3mm]

Notice that  Theorem~D can be viewed as a  strong topological version of the double ergodicity property of Poisson boundaries due to Burger-Monod (\cite{MR1911660}) and Kaimanovich (\cite{MR2006560}).

We next fix a base point $\xo\in\XX$, $\theta\in [0,\pi/2]$  and consider the cardinality of  the sets
$$ N_\theta^\eps(n):=
\{ \gamma\in\Gamma\;:\, n-1< d(\xo,\gamma\xo)\le n\,,\ \Big| \frac{d_2(p_2(\gamma \xo),p_2(\xo))}{d_1(p_1(\gamma\xo ),p_1(\xo))}-\tan\theta\Big|<\eps\}\,,$$
where  $\eps>0$ and $n\in\NN$ is large. This number counts orbit points with a correlation of distances to the origin in each factor given by approximately $\tan\theta$. We further define
$$\delta_\theta^\eps:=\lim\sup_{n\to\infty}\frac{\log N_\theta^\eps(n)}{n}\,, \quad \an\qquad \delta_\theta(\Gamma):=\liminf_{\eps\to 0} \delta_\theta^\eps\,.$$
$\delta_\theta(\Gamma)$ can be thought of as a function of $\theta\in [0,\pi/2]$ which describes the exponential growth rate of orbit points converging to limit points of slope $\theta$. It is an invariant of $\Gamma$ which carries more information than the critical exponent $\delta(\Gamma)$: the critical exponent is simply the maximum of $\delta_\theta(\Gamma)$ in $[0,\pi/2]$. 
As in  \cite{MR1933790} it will be convenient to study the homogeneous function 
$$\Psi_\Gamma:\RR_{\ge 0}^2 
\to \RR\,, \ H=(H_1,H_2)\mapsto \Vert H\Vert\cdot \delta_{\arctan(H_1/H_2)}(\Gamma)\,.$$ 
Similar to the case of symmetric spaces or Euclidean buildings of higher rank, we have the following: \\[2mm]
{\bf Theorem E}$\quad$ {\sl $\Psi_\Gamma$ is upper semi-continuous and concave.}\\[3mm]
One of the main applications of Theorem E is that it allows to construct generalized conformal densities on each $\Gamma$-invariant subset of the limit set as in \cite{MR2062761} and \cite{MR1935549} for higher rank symmetric spaces and Euclidean buildings. In a future work we will carry out this construction and relate $\delta_\theta(\Gamma)$ to the Hausdorff dimension of the limit set.

The paper is organized as follows: Section~2 recalls 
basic facts about Hadamard spaces  and rank one isometries. In Section~3  we collect  properties of products of Hadamard spaces. In Section~4 we study the structure of the limit set and prove Theorems A and B. 
Section~5 deals with properties of the set of regular axial isometries and contains the proofs of Theorems C and D.  In 
Section~6 we introduce and study the exponent of growth of slope $\theta$ for $\Gamma$. Finally, in Section~7 we construct a so-called generic product for $\Gamma$ in order to show that the function $\Psi_\Gamma$ 
is concave, and give the proof of Theorem E. \\[-2mm]

{\bf Acknowledgements:}  This paper was written during the author's stay at ETH Zurich. She warmly thanks Marc Burger and Alessandra Iozzi for inviting her, and  the FIM for its hospitality and the inspiring atmosphere. She is grateful to Pierre-Emmanuel Caprace for many helpful remarks and discussions, and to Fran{\c c}oise Dal'bo for her valuable comments on a first draft of the paper. She also thanks both 
referees for their useful suggestions and in particular a considerable simplification of the proof of Lemma~\ref{prodtopology}.

\section{Preliminaries}\label{Prelim}

The purpose of this section is to introduce some terminology and notation and to summarize basic results about Hadamard spaces and rank one isometries. 
The main references here  are \cite{MR1744486} and  \cite{MR1377265} (see also \cite{MR1383216}, and \cite{MR823981},\cite{MR656659} in the case of Hadamard manifolds). 

Let $(\XX,d)$ be a metric space. A {\hl geodesic path} joining $x\in\XX$ to $y\in\XX$  is a map $\sigma$ from a closed interval $[0,l]\subset \RR$ to $\XX$ \st $\sigma(0)=x$, $\sigma(l)=y$ and $d(\sigma(t), \sigma(t'))=|t-t'|$ for all $t,t'\in [0,l]$.  We will denote such a geodesic path $\sigma_{x,y}$. $\XX$ is called {\hl geodesic}, if any two points in $\XX$ can be connected by a geodesic path, if this path is unique, we say that $\XX$ is {\hl uniquely geodesic}. In this text $\XX$ will be a Hadamard space, i.e. a complete geodesic metric space in which all triangles satisfy the CAT$(0)$-inequality. This implies in particular that $\XX$ is simply connected and uniquely geodesic.   A {\hl geodesic} or {\hl geodesic line} in $\XX$ is a map $\sigma:\RR\to\XX$ \st $d(\sigma(t), \sigma(t'))=|t-t'|$ for all $t,t'\in\RR$, a {\hl geodesic ray} is a map $\sigma:[0,\infty)\to \XX$ \st $d(\sigma(t), \sigma(t'))=|t-t'|$ for all $t,t'\in [0,\infty)$. Notice that in the non-Riemannian setting completeness of $\XX$ does not imply that every geodesic path or ray can be extended to a geodesic, i.e. $\XX$ need not be geodesically complete.

From here on we will assume that $\XX$ is a locally compact 
Hadamard space.  The geometric boundary $\rand$ of
$\XX$ is the set of equivalence classes of asymptotic geodesic
rays endowed with the cone topology (see e.g. \cite[chapter~II]{MR1377265}). The action of the isometry group $\is(\XX)$ on $\XX$ naturally extends to an action by homeomorphisms on the geometric boundary. Moreover, since $\XX$ is locally compact, this boundary $\rand$ is compact and the space $\XX$ is a dense and open subset of the compact space $\ganz:=\XX\cup\rand$.  For $x\in\XX$ and $
\xi\in\rand$ arbitrary, there exists a  geodesic ray emanating from $x$ which belongs to the class of $\xi$. We will denote such a ray $\sigma_{x,\xi}$.

We say that two points $\xi$, $\eta\in\rand$ can be joined by a geodesic if there exists a geodesic $\sigma:\RR\to\XX$ \st $\sigma(-\infty)=\xi$ and $\sigma(\infty)=\eta$. It is well-known that if $\XX$ is  CAT$(-1)$, i.e. of negative Alexandrov curvature bounded above by $-1$, then every pair of distinct points in the geometric boundary can be joined by a geodesic. This is not true in general. For convenience we therefore define the {\hl visibility set at infinity} $\horinf(\xi)$ of a point  $\xi\in\rand$ as  the set of points in the geometric boundary which can be joined to $\xi$ by a geodesic, i.e.
\begin{equation}\label{visinf}
\horinf(\xi):=\{\eta\in\rand\;|\ \exists\ \mbox{geodesic}\ \sigma\  \st 
\sigma(-\infty)=\xi\,,\,\sigma(\infty)=\eta\}\,.
\end{equation}

Let $x, y\in \XX$, $\xi\in\rand$ and $\sigma$ a geodesic ray in the
class of $\xi$. We put 
\begin{equation}\label{buseman}
 \bs_{\xi}(x, y)\,:= \lim_{s\to\infty}\big(
d(x,\sigma(s))-d(y,\sigma(s))\big)\,.
\end{equation}
This number is independent of the chosen ray $\sigma$, and the
function
\be \bs_{\xi}(\cdot , y):
\quad \XX &\to & \RR\\
x &\mapsto & \bs_{\xi}(x, y)\ee
is called the {\hl Busemann function} centered at $\xi$ based at $y$ (see also \cite{MR1377265}, chapter~II). For any $x,y,z\in\XX$, $\xi\in\rand$ and $g\in\is(\XX)$ the Busemann function satisfies 
\begin{eqnarray}
|\bs_{\xi}(x, y)|&\le &d(x,y)\label{bounded}\\
\bs_{\xi}(x, z)&=&\bs_{\xi}(x, y)+\bs_{\xi}(y,z)\label{cocycle}\\
\bs_{g\cdot\xi}(g\at x,g\at y) & =& \bs_{\xi}(x, y)\,.\nonumber
\end{eqnarray}
A geodesic $\sigma: \RR\to\XX$  is said to {\hl bound a flat half-plane}  if there exists a closed convex subset $i([0,\infty)\times \RR)$ in $\XX$ isometric to $[0, \infty)\times\RR$ \st $\sigma(t)=i(0,t)$ for all $t\in \RR$. Similarly, a geodesic $\sigma: \RR\to\XX$  bounds a {\hl flat strip} of width $c>0$   if there exists a closed convex subset $i([0,c]\times \RR)$ in $\XX$ isometric to $[0, c]\times\RR$ \st $\sigma(t)=i(0,t)$ for all $t\in \RR$. 
We call a geodesic $\sigma:\RR\to\XX$ a {\hl rank one geodesic} if  $\sigma$ does not bound a
flat half-plane.  

The following important lemma states that even though we cannot join any two distinct points in the geometric boundary of $\XX$, given a rank one geodesic we can at least join points in a neighborhood of its extremities. More precisely, we have the following well-known
\begin{lem}\label{joinrankone} (\cite{MR1377265}, Lemma III.3.1)\ 
Let $\sigma:\RR\to\XX$ be a rank one geodesic which does not bound a flat strip of width $c$. Then there are neighborhoods $U$ of $\sigma(-\infty)$ and $V$ of $\sigma(\infty)$ in $\ganz$ \st for any $\xi\in U$ and $\eta \in V$ there exists a rank one geodesic joining $\xi$ and $\eta$. For any such geodesic $\sigma'$ we have $d(\sigma', \sigma(0))\le c$. 
\end{lem}
Moreover, we will need the following technical lemma which immediately follows from Lemma~4.3 and Lemma~4.4 in \cite{MR1383216}.
\begin{lem}\label{convergence}
Let $\sigma:\RR\to\XX$ be a rank one geodesic and put $y:=\sigma(0)$, $\eta:=\sigma(\infty)$. Then for any $T\gg1$, $\eps>0$  there exists a neighborhood $U$ of $\sigma(-\infty)$ in $\ganz$ and a number $R>0$ \st for any $x\in \XX$ with $d(x,\sigma)>R$ or $x\in U$ we have
$$ d(\sigma_{x,y}(t),\sigma_{x,\eta}(t))\le \eps\ \quad\mbox{for all}\ \ t\in [0,T]\,.$$
\end{lem}

The following kind of isometries will play a central role in the sequel. 
\begin{df} \label{axialisos}
An isometry $h$ of $\XX$ is called {\hd axial}, if there exists a constant
$l=l(h)>0$ and a geodesic $\sigma$ \st
$h(\sigma(t))=\sigma(t+l)$ for all $t\in\RR$. We call
$l(h)$ the {\hd translation length} of $h$, and $\sigma$
an {\hd axis} of $h$. The boundary point
$h^+:=\sigma(\infty)$ is called the {\hd attractive fixed
point}, and $h^-:=\sigma(-\infty)$ the {\hd repulsive fixed
point} of $h$. We further put
$\Ax(h):=\{ x\in\XX\;|\, d(x,h x)=l(h)\}$.
\end{df}
We remark that $\Ax(h)$ consists of the union of parallel geodesics
translated by $h$, and 
$\overline{\Ax(h)}\cap\rand$ is exactly the set of fixed points of
$h$. Moreover, we have the following easy formula for the translation length of an axial isometry in terms of Busemann functions.
\begin{lem}\label{transbus}
If $h$ is an axial isometry with attractive and repulsive fixed points $h^+$, $h^-$ then its translation length is given by
$$l(h)=\bs_{h^+}(x,h  x)=\bs_{h^-}(h  x,x)\,,\quad \mbox{where $x\in\XX$ is arbitrary}\,.$$
\end{lem}
\prf\ Let $x$, $y\in\XX$ arbitrary. Then by the cocycle identity~(\ref{cocycle}) and the fact that $h$ fixes $h^+$ and $h^-$ 
\be \bs_{h^+}(x,h  x) & = &\bs_{h^+}(x,y)+\bs_{h^+}(y,h  y)+\bs_{h^+}(h  y ,h  x)\\
&=& \bs_{h^+}(x,y)+\bs_{h^+}(y,h  y)+\underbrace{\bs_{h^+}(y, x)}_{=-\bs_{h^+}(x,y)}=\bs_{h^+}(y,h  y)\,,
\ee
and similarly $\bs_{h^-}(h  x,x)=\bs_{h^-}(h  y,y)$. So the terms on the right-hand side are independent of $x\in\XX$, and  choosing $x\in\Ax(h)$ yields the claim.\qed\\[3mm]
Following the definition in \cite{GR/0702274} and \cite{GR/0809.0470} we will call two axial isometries $g$, $h\in\is(\XX)$ {\hl independent} if for any given $x\in\XX$ the map 
$$ \ZZ\times\ZZ\to [0,\infty)\,,\ (m,n)\mapsto d(g^m   x, h^m  x)$$
is proper.
\begin{df}
An axial isometry  is called {\hd rank one} if it possesses a rank one axis. 
\end{df}
Notice that if $h$ is rank one, then $h^+$ and $h^-$ are the only fixed points of $h$. Moreover, it is easy to verify that two rank one elements $g$, $h\in\is(\XX)$ are independent if and only if $\{g^+,g^-\}\cap\{h^+,h^-\}=\emptyset$. 
Let us recall some properties of rank one isometries.
\begin{lem}\label{dynrankone}(\cite{MR1377265}, Lemma III.3.3)\ 
Let $h$ be a rank one isometry. Then
\begin{enumerate}
\item[(a)] $\horinf(h^+)=\rand\setminus\{h^+\}$, 
\item[(b)] any geodesic joining a point $\xi\in\rand\setminus\{h^+\}$ to $h^+$ is rank one, 
\item[(c)] given neighborhoods $U$ of $h^-$ and $V$ of $h^+$ in $\ganz$ 
there exists $N_0\in\NN$ \st\\
 $h^{-n}(\ganz\setminus V)\subset U$ and
$h^{n}(\ganz\setminus U)\subset V$ for all $n\ge N_0$.
\end{enumerate}
\end{lem}
In particular, by (c) we have $\lim_{n\to\infty} h^{-n}\xi=h^-$ for any $\xi\in\horinf(h^+)$.
The following lemma will be central for the proof of Theorem~\ref{propertiesofPGamma}. 
\begin{lem}\label{elementsarerankone} (\cite{MR1377265}, Lemma III.3.2)\
Let $\sigma:\RR\to\XX$ be a rank one geodesic, and $(\gamma_n)\subset\is(\XX)$ a sequence of isometries \st $\gamma_n x\to \sigma(\infty)$ and $\gamma_n^{-1}x\to\sigma(-\infty)$ for one (and hence any) $x\in\XX$. Then for $n$ sufficiently large, $\gamma_n$ is axial and possesses an axis $\sigma_n$ \st $\sigma_n(\infty)\to\sigma(\infty)$ and $\sigma_n(-\infty)\to\sigma(-\infty)$.
\end{lem}
The following proposition is a generalization of Lemma~4.1 in \cite{MR1703039}. It gives a relation between the geometric length and the combinatorial length of words in a free group on two generators which will be a clue to the proof of Theorem~\ref{propertiesofPGamma}. Our proof here involves a new idea since 
F.~Dal'bo's proof 
is based  on  the fact that  $\XX$ is CAT$(-1)$ and hence triangles in $\XX$ are thinner than the corresponding triangles in hyperbolic space. If $g,h$ generate a free group we say that a word $\gamma=s_1^{k_1}s_2^{k_2}\cdots s_n^{k_n}$ with $s_i\in \{g,g^{-1}, h,h^{-1}\}$ and $k_i\in\NN\setminus\{0\}$, $i\in\{1,2,\ldots n\}$   is {\hl cyclically reduced} if $s_{i+1}\notin \{s_i, s_i^{-1}\}$, $i\in\{1,2,\ldots  n-1\}$, and $s_n\ne s_1^{-1}$.  
\begin{prp}\label{combgeomlength}
Suppose  $g$ and $h$ are rank one elements in $\is(\XX)$ with pairwise distinct fixed points. Then there exists $N\in\NN$ and $C>0$ \st %the group generated by $g^N$ and $h^N$ is free, and 
for all $n\in\NN$ and any cyclically reduced word $\gamma=s_1^{k_1}s_2^{k_2}\cdots s_n^{k_n}$ with $s_i\in S:=\{g^N,g^{-N}, h^N,h^{-N}\}$ and $k_i\in\NN\setminus\{0\}$, $i\in\{1,2,\ldots n\}$, we have
$$ \big|l(\gamma)-\sum_{i=1}^n k_i l(s_i)\big|\le C\cdot n\,.$$
\end{prp}
\prf\  We fix some base point $\xo\in\XX$. For $\eta\in \{g^-,g^{+},h^-,h^{+}\}$ let $U(\eta)\subset \ganz$ be a small neighborhood of $\eta$ with $\xo\notin U(\eta)$ \st all $U(\eta)$ are pairwise disjoint, and $c>0$ a constant \st any pair of points  in distinct neighborhoods can be joined by a rank one geodesic $\sigma'$ with $d(\xo,\sigma')\le c$. This is possible by Lemma~\ref{joinrankone}. According to Lemma~\ref{elementsarerankone} there exist neighborhoods $W(\eta)\subseteq U(\eta)$, $\eta\in \{g^-,g^{+},h^-,h^{+}\}$, \st every $\gamma\in\Gamma$ with $\gamma\xo\in W(\eta)$, $\gamma^{-1}\xo\in W(\zeta)$, $\zeta\neq \eta$, is rank one with $\gamma^+\in U(\eta)$ and $\gamma^-\in U(\zeta)$. Moreover, by Lemma~\ref{dynrankone} (c) there exists $N\in\NN$ \st for all $\gamma\in\{g,g^{-1},h,h^{-1}\}$
\begin{equation}\label{pingpong} 
\gamma^N\big(\ganz\setminus W(\gamma^{-})\big)\subseteq W(\gamma^+)\,.
\end{equation}
We put $S:=\{g^N,g^{-N},h^N,h^{-N}\}$ and consider a cyclically reduced word $
\gamma=s_1^{k_1}s_2^{k_2}\cdots s_n^{k_n}$ with $s_i\in S$ and $k_i\in\NN\setminus\{0\}$, $i\in\{1,2,\ldots n\}$. By the choice of $N$ and (\ref{pingpong}) we have $\gamma\xo\in W(s_1^+)$ and $\gamma^{-1}\xo\in W(s_n^{-})\neq W(s_1^+)$ since $s_1\neq s_n^{-1}$. Therefore $\gamma$ is rank one with $\gamma^+\in U(s_1^+)$ and $\gamma^-\in U(s_n^{-})$. Choosing a point $x\in\Ax(\gamma)$ with $d(\xo,x)\le c$ we get
\begin{equation}\label{translengthdistance}
 l(\gamma)\le d(\xo,\gamma\xo)\le d(\xo,x)+d(x,\gamma x)+ d(\gamma x,\gamma\xo)\le l(\gamma)+2c\,.
\end{equation}
Similarly $l(s_i^{k_i})\le d(\xo,s_i^{k_i}\xo)\le l(s_i^{k_i})+2c$ for $i\in\{1,2,\ldots n\}$.

For $i\in\{1,2,\ldots n\}$ we abbreviate $\gamma_i:=s_i^{k_i}s_{i+1}^{k_{i+1}}\cdots s_n^{k_n}$. Then $\gamma_2\xo\in W(s_2^+)$, $s_1^{-k_1}\xo\in W(s_1^{-})\neq W(s_2^+)$, so there exists a geodesic $\sigma_2$ joining $\gamma_2\xo$ to $s_1^{-k_1}\xo$ with $d(\xo,\sigma_2)\le c$. If $y$ denotes a point on $\sigma_2$ with $d(\xo,y)\le c$ we obtain
$$ d(s_1^{k_1}s_2^{k_2}\cdots s_n^{k_n}\xo,\xo)= d(\gamma_2\xo, s_1^{-k_1}\xo)\le d(\gamma_2\xo, y)+d(y,s_1^{-k_1}\xo)$$ 
which proves $|d(\gamma\xo,\xo)-d(\xo,s_1^{k_1}\xo)-d(\xo,\gamma_2\xo)|\le 2c$. Applying the same arguments to $\gamma_i$ for $i\ge 2$ and using the fact that $s_{i+1}\ne s_i^{-1}$ we deduce $|d(\gamma_i\xo,\xo)-d(\xo,s_i^{k_i}\xo)-d(\xo,\gamma_{i+1}\xo)|\le 2c$. Therefore 
$$ \big| d(\xo,\gamma\xo)-\sum_{i=1}^n d(\xo,s_i^{k_i}\xo)\big|\le 2(n-1)c $$
and, using (\ref{translengthdistance}),  we conclude
$ \big|l(\gamma)-\sum_{i=1}^n k_i l(s_i)\big|\le 4c\cdot n$. It remains to set $C:=4c$. \qed\\[3mm]
Moreover, the following generalization of Lemma~1.4 (2) in \cite{DalboKim} will also be needed in the proof of Theorem~\ref{propertiesofPGamma}:
\begin{lem}\label{combgeomlen}
Suppose  $g$ and $h$ are rank one elements in $\is(\XX)$ with $g^+=h^+$. Then there exists $N\in\NN$  \st for all $n,m \in\NN\setminus\{0\}$  the isometry $g^{Nn}h^{Nm}$ is rank one and $$l(h^{Nn}g^{Nm})=Nn\;l(h)+Nm\; l(g)\,.$$
\end{lem}
\prf\ As in the proof of the previous proposition we fix some base point $\xo\in\XX$ and let  $U(\eta)\subset \ganz$  be a  small neighborhood of $\eta\in \{g^-,g^{+},h^-\}$ with $\xo\notin U(\eta)$  \st all $U(\eta)$ are pairwise disjoint. Notice that by our assumption we may set $U(h^+):=U(g^+)$. 
Fix  neighborhoods $W(\eta)\subseteq U(\eta)$, $\eta\in \{g^-,g^{+},h^-, h^+\}$, \st every $\gamma\in\Gamma$ with $\gamma\xo\in W(\eta)$, $\gamma^{-1}\xo\in W(\zeta)$, $\zeta\neq \eta$, is rank one with $\gamma^+\in U(\eta)$ and $\gamma^-\in U(\zeta)$, and   $N\in\NN$ \st for all $\gamma\in\{g,g^{-1},h,h^{-1}\}$
$$%\begin{equation}
\gamma^N\big(\ganz\setminus W(\gamma^{-})\big)\subseteq W(\gamma^+)\,.
$$%\end{equation}
Then for $n,m\in\NN\setminus\{0\}$ $h^{Nn} g^{Nm}\xo\in W(h^+)$ and $(h^{Nn}g^{Nm})^{-1}\xo\in W(g^-)\ne W(h^+)$, hence $\gamma:=h^{Nn} g^{Nm}$ is rank one with $\gamma^+\in U(h^+)$ and $\gamma^-\in U(g^-)$. Furthermore, $\gamma h^+= h^{Nn}g^{Nm} h^+=h^+$ implies that $h^+$ is one of the two fixed points of $\gamma$, hence $\gamma^+=h^+=g^+$. We conclude using Lemma~\ref{transbus} and the cocycle identity (\ref{cocycle})
\be
l(\gamma)&=& \bs_{\gamma^+}(\xo,\gamma\xo) =\bs_{h^+}(\xo,h^{Nn}\xo)+\bs_{h^+}(h^{Nn}\xo, h^{Nn}g^{Nm}\xo)\\
&=& l(h^{Nn})+\bs_{g^+}(\xo,g^{Nm}\xo)= Nn\;l(h)+Nm\;l(g)\,.\hspace{4.7cm} \Box\ee

If  $\Gamma$ is a group acting by isometries on a locally compact Hadamard space $\XX$ we define its 
{\hl geometric limit set}  by
$\Lim:=\overline{\Gamma\at x}\cap\rand$, where $x\in \XX$ is arbitrary.

From here on we let $\Gamma\subset\is(\XX)$ be a (not necessarily discrete) group which possesses a rank one element $h$. Denote $\sigma$ an axis of $h$ and put $\xo:=\sigma(0)$. The idea of proof of the following three facts is due to W.~Ballmann (see e.g. the proof of Theorem~2.8 in \cite{MR656659}). We include complete proofs for the convenience of the reader. 
\begin{lem}\label{goesinnbh}
If $\Gamma$ does not globally fix a point in $\rand$, then for any neighborhood $V$ of $\xi\in\Lim$ in $\ganz$ there exists $\gamma\in\Gamma$ \st $\gamma h^+\in V$.
\end{lem}
\prf\ Choose $(\gamma_n)\subset\Gamma$ \st $\gamma_n\xo\to\xi$ as $n\to\infty$. Passing to a subsequence if necessary we may assume that $\gamma_n^{-1}\xo$ converges to a point $\zeta\in\Lim$ as $n\to\infty$. Let $T\gg1$ and $\eps>0$ be arbitrary. By Lemma~\ref{convergence} there exist a constant $R>0$ and a neighborhood $U$ of $h^-=\sigma(-\infty)$ in $\ganz$ \st for any $x\in\XX$ with $d(x,\sigma)>R$ or $x\in U$ we have
$ d(\sigma_{x,\xo}(t),\sigma_{x,h^+}(t))\le \eps/2$ for all $t\in [0,T]$\,.

We first treat the case $\zeta\notin\{h^+,h^-\}$. Then for $n$ sufficiently large we have $d(\gamma_n^{-1}\xo,\sigma)>R$ and $d(\sigma_{\xo,\gamma_n\xo}(t),\sigma_{\xo,\xi}(t))\le \eps/2$ for $0\le t\le T$. We conclude that for $t\in [0,T]$
\be d(\sigma_{\xo,\xi}(t),\sigma_{\xo,\gamma_n h^+}(t)) &\le & d(\sigma_{\xo,\xi}(t),\sigma_{\xo,\gamma_n\xo}(t))+d(\sigma_{\xo,\gamma_n\xo}(t),\sigma_{\xo,\gamma_n h^+}(t))\\ 
&\le & \frac{\eps}2 +d(\gamma_n \sigma_{\gamma_n^{-1}\xo,\xo}(t),\gamma_n \sigma_{\gamma_n^{-1}\xo,h^+}(t))\\
&=& \frac{\eps}2+\underbrace{d(\sigma_{\gamma_n^{-1}\xo,\xo}(t),\sigma_{\gamma_n^{-1}\xo,h^+}(t))}_{\le \eps/2}\le \eps\,,
\ee
which proves the assertion in this case.

If $\zeta=h^-$ then $\gamma_n^{-1}\xo\in U$ for $n$ sufficiently large, hence by Lemma~\ref{convergence}  and the above inequalities the claim also holds.

 Now assume that $\zeta= h^+$. Since $\Gamma$ does not globally fix a point in $\rand$ there exists $\varphi\in\Gamma$ \st $\varphi\zeta\ne h^+$. Then, replacing $\gamma_n$ by $\gamma_n \varphi^{-1}$ and using the fact that $\gamma_n\varphi^{-1}\xo\to \xi$ and $\varphi\gamma_n^{-1}\xo\to \varphi\zeta\ne h^+$, we are in one of the cases above. Hence the assertion follows. \qed\\[3mm]
The following result  will be one of the key lemmas for the product case in Section~\ref{limsetdiscreteisom}.
\begin{lem}\label{movepoints}
If $\Gamma$ does not globally fix a point in $\rand$ and $\#\Lim=\infty$, then for all $\xi$, $\eta$, $\zeta\in\Lim$ there exists $\gamma\in\Gamma$ \st $\gamma\xi\neq\zeta$ and $\gamma\xi\neq\eta$.
\end{lem} 
\prf\  If $\xi\in\Lim\setminus\{\zeta,\eta\}$ we can take $\gamma=e$ (the identity in $\Gamma$). 

Suppose now $\xi=\zeta\neq\eta$ and $\xi\notin\{h^+,h^-\}$. Then $h^n\xi\to h^+$ as $n\to\infty$. If $\eta\neq h^+$, let $V$ be a neighborhood of $h^+$ disjoint from $\xi$, $\eta$. Then there exists $N\in\NN$ \st $h^n\xi\in V$ for all $n\ge N$, in particular $h^N\xi\neq \xi=\zeta$ and $h^N\xi\neq \eta$. If $\eta=h^+$ we choose a neighborhood $V$ of $h^+$ disjoint from $\xi$ and let $N\in\NN$ \st $h^n\xi\in V$ for all $n\ge N$. If  $h^n\xi=\eta=h^+$ for all $n\ge N$, then $\xi$ is a fixed point of $h$ which is a contradiction to $\xi\notin\{h^+,h^-\}$. Hence there exists $n\ge N$ \st $h^n\xi\ne\eta$ and $h^n\xi\ne\xi=\zeta$.

If $\xi=\zeta=h^+$, $\eta\ne h^+$ we choose a point in $\Lim\setminus\{h^+,\eta\}$ and a neighborhood $V$ of this point disjoint from $\{h^+,\eta\}$. By Lemma~\ref{goesinnbh} there exists $\gamma\in\Gamma$ \st $\gamma h^+\in V$, in particular $\gamma h^+\notin \{h^+,\eta\}$. 

Replacing $h$ by $h^{-1}$ in the previous argument yields the assertion for
the case $\xi=\zeta=h^-$, $\eta\ne h^-$. 

By symmetry, the claim also holds for $\xi=\eta\neq\zeta$.

The remaining case is $\xi=\zeta=\eta$. Since $\Gamma$ does not globally fix a point in $\rand$, there exists $\gamma\in\Gamma$ \st $\gamma\xi\neq\xi$. \qed \\[3mm]
In the case of discrete groups, the following result is part of Theorem~2.8 in \cite{MR656659}.  Since we are dealing here with possibly non-discrete groups we have to add the condition that $\Gamma$ does not globally fix a point in $\rand$. This excludes for example the case of a group consisting of infinitely many rank one elements with a common fixed point at infinity.
\begin{prp}\label{minimalityinfactor}
If $\Gamma$ does not globally fix a point in $\rand$ and $\#\Lim=\infty$, then the limit set $\Lim$ is minimal, i.e. the smallest non-empty $\Gamma$-invariant closed subset of $\rand$. 
\end{prp}
\prf\ We first notice that every non-empty $\Gamma$-invariant closed subset $A$ of $\rand$ contains a limit point: Indeed, if $\xi\in A$, then either $\xi=h^+$ or $\xi\in\horinf(h^+)$. So either $A$ 
contains the limit point $h^+$ or the point $h^-=\lim_{n\to\infty} h^{-n}\xi$. 

Next we fix $\xi\in\Lim$ and let $\eta\in\Lim$ be arbitrary. Our goal is to show that $\eta\in\overline{\Gamma\cdot\xi}$.

Let $U\subset\rand$ be an arbitrary neighborhood of $\eta$. By Lemma~\ref{goesinnbh} there exists $\gamma\in\Gamma$ \st $\gamma h^+\in U$. Hence if $\xi\neq\gamma h^-$ we have by the dynamics of rank one isometries Lemma~\ref{dynrankone} (c) $(\gamma h\gamma^{-1})^n\xi\in U$ for $n$ sufficiently large. If $\xi=\gamma h^-$, there exists $\varphi\in\Gamma$ \st $\varphi\xi\neq \gamma h^-$ by Lemma~\ref{movepoints}. Then $(\gamma h\gamma^{-1})^n\varphi\xi\in U$ for $n$ sufficiently large. \qed

\section{Products of Hadamard spaces}\label{prodHadspaces}

Now let $(\XX_1,d_1)$, $(\XX_2,d_2)$ be locally compact  
Hadamard spaces,  and $\XX=\XX_1\times \XX_2$  the product space endowed with the product distance $d=\sqrt{d_1^2+d_2^2}$. Notice that such a product is again a locally compact 
Hadamard space.  To any pair of points $x=(x_1,x_2)$, $z=(z_1,z_2)\in \XX$ we associate the vector 
$$H(x,z):= \left(\begin{array}{c} d_1(x_1,z_1)\\ d_2(x_2,z_2)\end{array}\right)\in \RR^2\,,$$ 
which we call the {\hl distance vector} of the pair $(x,z)$. If $z\neq x$ we further define the {\hl direction} of $z$  \wrt $x$ by 
$$\theta(x,z):=\arctan\frac{d_2(x_2,z_2)}{d_1(x_1,z_1)}\,.$$
Notice that we have $$H(x,z) = d(x, z)\left(\begin{array}{c} \cos \theta(x,z)\\ \sin\theta(x,z)\end{array}\right)\,,$$ in particular $\Vert H(x,z)\Vert =d(x,z)$, where $\Vert \cdot\Vert$ denotes the Euclidean norm in $\RR^2$. 

Denote $p_i:\XX\to \XX_i$, $i=1,2$, the natural projections. Every 
geodesic path $\sigma:[0,l]\to\XX\,$ can be written as a product $\sigma(t)=(\sigma_1(t \cos \theta), \sigma_2(t \sin\theta))$, where $\theta\in [0,\pi/2]$ and $\sigma_1:[0,l \cos\theta]\to\XX_1$, $\sigma_2:[0,l \sin\theta]\to\XX_2$ are geodesic paths in $\XX_1$, $\XX_2$. $\theta$ equals the direction of $\sigma(l)$ with respect to $\sigma(0)$ and  is called the {\hl slope of $\sigma$}.  We say that a geodesic path $\sigma$ is  {\hl regular} if its slope is contained in the open interval $(0,\pi/2)$. In other words, $\sigma$ is regular if neither $p_1(\sigma([0,l]))$ nor $p_2(\sigma([0,l]))$ is a point.

If  $x\in\XX$  and $\sigma:[0,\infty)\to\XX$ is an arbitrary geodesic ray, then by elementary geometric estimates one has  the relation 
\begin{equation}\label{RelSlopeDir}
 \theta = \lim_{t\to\infty} \theta(x,\sigma(t)) 
 \end{equation}
 between the slope $\theta$ of $\sigma$ and the directions of $\sigma(t)$, $t>0$, \wrt $x$. 
 Similarly, one can easily show that any two geodesic rays representing the same (possibly singular) point in the geometric boundary necessarily have the same slope. So  we may define the slope $\theta(\tilde\xi)$ of a point $\tilde\xi\in\rand$ as the slope of an arbitrary geodesic ray representing $\tilde\xi$. 
 
 Moreover, two regular geodesic rays  $\sigma$, $\sigma'$  with the same slope represent the same point in the geometric boundary if and only if  $\sigma_1(\infty)=\sigma_1'(\infty)$ and  $\sigma_2(\infty)=\sigma_2'(\infty)$.  The regular geometric boundary $\regrand$ of $\XX$ is defined as the set of equivalence classes of regular geodesic rays and hence  is  homeomorphic  to $\rand_1\times \rand_2\times (0,\pi/2)$. 

 If $\gamma\in\is(\XX_1)\times \is(\XX_2)$ , then 
 the slope of $\gamma\at\tilde\xi$ equals the slope of $\tilde\xi$. In other words, if $\rand_\theta$ denotes the set of points in the geometric boundary of slope $\theta\in [0,\pi/2]$, then $\rand_\theta$ is invariant by the action of $\is(\XX_1)\times\is(\XX_2)$. Notice that points in $\singrand:= (\rand)_0\cup(\rand)_{\pi/2}$ are equivalence classes of geodesic rays which project to a point in one of the factors of $\XX$. Hence $(\rand)_0$ is homeomorphic to $\partial \XX_1$ and $(\rand)_{\pi/2}$ is homeomorphic to $\rand_2$. 
If  $\theta\in (0,\pi/2)$, then the set $\rand_\theta\subset\regrand$ is homeomorphic to the product $\rand_1\times\rand_2$. 

In the case of symmetric spaces and Bruhat-Tits buildings of higher rank there is a well-known notion of Furstenberg boundary, which -- for a product of rank one spaces -- coincides with the product of the geometric boundaries. In our more general setting we therefore choose to call  the product $\rand_1\times \rand_2$ endowed with the product topology the {\hl Furstenberg boundary} $\Frand$ of $\XX$.
Using the above parametrization of $\regrand$ we have a natural projection
$$ \hspace{1cm} \ba{rcl}\pi^F\,:\qquad\ \regrand&\to& \Frand\\
(\xi_1,\xi_2,\theta)&\mapsto & (\xi_1,\xi_2)\, \ea $$
and a natural action of the  group $\is(\XX_1)\times\is(\XX_2)$ by homeomorphisms on the Furstenberg boundary of $\XX=\XX_1\times\XX_2$.

We have the following important  lemma concerning the topology of $\ganz$. Although elementary, we include the proof for the convenience of the reader.
\begin{lem}\label{prodtopology}
Suppose $(y_n)\subset\XX$ is a sequence converging to a point $\tilde\eta\in\rand_\theta$ for some $\theta\in [0,\pi/2]$. Then for any $x\in\XX\,$ we have $\theta(x,y_n)\to \theta$ as $n\to\infty$. 
\end{lem}
\prf\ First notice that if $\sigma$ is a geodesic emanating from $x$, then $\theta(x,\sigma(t))$ 
does not depend on $t$. We define $\sigma$ as a geodesic ray joining $x$ to
$\tilde\eta$, so in particular $\sigma$ has slope $\theta$ and $\theta(x,\sigma(t))=\theta$ for all $t>0$. Without loss of generality we may assume that $d(x,y_n)\ge 1$ for all $n\in\NN$. It therefore remains to prove that $\theta(x,\sigma_{x,y_n}(1))=\theta(x,y_n)$ converges to $\theta(x,\sigma(1))=\theta$  as $n$ tends to infinity. This is clear since $\sigma_{x,y_n}(1)$  converges to $\sigma(1)$  and since the
map $z\mapsto \theta(x,z)$  is continuous on every sphere around $x$.\qed\\

Recall the definition of visibility set at infinity  $\horinf(\tilde\xi)$  of a point $\tilde\xi\in\rand$ from~(\ref{visinf}). It is easy to see that a point $\tilde\eta\in\rand$ cannot belong to $\horinf(\tilde\xi)$ if the slope of $\tilde\eta$ is different from the slope of $\tilde\xi$. This motivates the following less restrictive definition for pairs of points in the Furstenberg boundary: We say that $\xi=(\xi_1,\xi_2)$ and $\eta=(\eta_1,\eta_2) \in\Frand$ are {\hl opposite}  if $\xi_1$ and $\eta_1$ can be joined by a geodesic in $\XX_1$, and $\xi_2$, $\eta_2$ can be joined by a geodesic in $\XX_2$. Moreover, 
the {\hl Furstenberg visibility set} $ \horF(\xi)$ of a point $\xi=(\xi_1,\xi_2)\in\Frand$  is defined as the set of points in $\Frand$ which are opposite to $\xi$, i.e. 
\begin{equation}\label{visFequiv}
\horF(\xi)=\{(\eta_1,\eta_2)\in\Frand:\eta_1\in\horinf(\xi_1)\ \an\ \  \eta_2\in\horinf(\xi_2)\}\,.
\end{equation}
In particular, for any $\tilde \xi\in\regrand$ with $\pi^F(\tilde \xi)=\xi$ one has $\horF(\xi)=\pi^F\big(\horinf(\tilde\xi)\big)$. So we may alternatively define the Furstenberg visibility set of a point $\xi\in\Frand$ via
\begin{equation}\label{visF}
 \horF(\xi):=\pi^F\big(\horinf(\tilde\xi)\big)\,,\quad \mbox{where}\ \  \tilde\xi\in (\pi^F)^{-1}(\xi)\ \mbox{is arbitrary}\,.
\end{equation}

Moreover, in the particular case that both  $\XX_1$, $\XX_2$ are CAT$(-1)$,  for $\tilde\xi=(\xi_1,\xi_2,\theta)\in\regrand$ we have
$$\horinf(\tilde\xi)=\{(\eta_1,\eta_2,\theta)\in\regrand : \eta_1\neq \xi_1 \mbox{ and}\ \eta_2\neq \xi_2\}\,,$$
and $(\xi_1,\xi_2)$, $(\eta_1,\eta_2)\in\Frand$ are opposite if and only if $\xi_1\neq \eta_1$ and $\xi_2\neq \eta_2$.

\section{The structure of the limit set}\label{limsetdiscreteisom}

Recall that the geometric limit set of a  group $\Gamma$ acting  by isometries on a locally compact 
Hadamard space is defined by
$\Lim:=\overline{\Gamma\at x}\cap\rand$, where $x\in \XX$ is arbitrary. In this section we will investigate the structure of the geometric limit set of certain groups $\Gamma\subset\is(\XX_1)\times\is(\XX_2)\subseteq\is(\XX)$ acting properly discontinuously on the product $\XX$ of two locally compact  
Hadamard spaces $\XX_1$, $\XX_2$. By abuse of notation we denote $p_i:\Gamma\to \is(\XX_i)$, $i=1,2$, the natural projections and put $\Gamma_i:=p_i(\Gamma)$, $i=1,2$. Notice that $\Gamma_i$ need not act properly discontinuously on $\XX_i$. As in \cite{DalboKim} for $i\in\{1,2\}$ we call $\Gamma_i$ {\hl strongly non-elementary} if it does not globally fix a point in $\rand_i$ and $\#L_{\Gamma_i}$ is infinite.

From here on we assume that $\Gamma\subset\is(\XX_1)\times\is(\XX_2)$ acts properly discontinuously, the projections $\Gamma_1$, $\Gamma_2$ are strongly non-elementary, 
and $\Gamma$ contains an isometry $h$ \st $h_1:= p_1(h)$, $h_2:=p_2(h)$ are rank one elements in $\Gamma_1$, $\Gamma_2$ respectively. Such an isometry $h$ will be called 
{\hl regular axial} and we will denote $\widetilde{h^+}$ its attractive fixed point in $\regrand$ and $h^+:=\pi^F(\widetilde{h^+})=(h_1^+, h_2^+)$. Notice that (\ref{visFequiv}) and Lemma~\ref{dynrankone} (a)  imply  
\begin{equation}\label{Fregaxial}
\horF(h^+)=\{(\xi_1,\xi_2)\in\Frand : \xi_1\ne h_1^+,\, \xi_2\ne h_2^+\}\,.
\end{equation}
Moreover, by Lemma~\ref{dynrankone} (c) we have $\lim_{n\to\infty} h^{-n}\xi=h^-$ for all $\xi\in\horF(h^+)$.

We remark that the existence of a regular axial element in $\Gamma$ imposes severe restrictions on the spaces $\XX_1$ and $\XX_2$. For example, neither $\XX_1$ nor $\XX_2$ can be a higher rank symmetric space 
or Euclidean building. However, as mentioned in the introduction the buildings associated to Kac-Moody groups over finite fields, Riemannian universal covers of geometric rank one manifolds  and CAT$(-1)$-spaces 
such as locally finite trees or manifolds of pinched negative curvature are natural 
examples of possible factors. 

For convenience we define the {\hl Furstenberg limit set} of $\Gamma$ by $\Flim:=\pi^F(\Lim\cap\regrand)$.  It is clearly a subset of the product $L_{\Gamma_1}\times L_{\Gamma_2}\subseteq\Frand$. Using our Lemma~\ref{movepoints} the proof of the following important lemma is as for Lemma~2.2 in \cite{DalboKim}.
\begin{lem}\label{goodposition}
For any $\xi=(\xi_1,\xi_2)$,  $\eta=(\eta_1,\eta_2)\in L_{\Gamma_1}\times L_{\Gamma_2}$ there exists $\gamma=(\gamma_1,\gamma_2)\in\Gamma$ \st $\gamma_1\xi_1\neq \eta_1$ and $\gamma_2\xi_2\neq\eta_2$. 
\end{lem}
\prf\  We first treat the case $\xi_1=\eta_1$ and $\xi_2\neq \eta_2$. Choose $\gamma=(\gamma_1,\gamma_2)$, $\varphi=(\varphi_1,\varphi_2)\in\Gamma$ \st $\gamma_1\xi_1\ne \eta_1$ and $\varphi_2\xi_2\notin\{\xi_2,\eta_2\}$. This is possible by Lemma~\ref{movepoints}. If $\gamma_2\xi_2\neq \eta_2$, $\gamma$ is the desired element, if $\varphi_1\xi_1\neq \eta_1$, then $\varphi$ is. 

Suppose now $\gamma_2\xi_2=\eta_2$ and $\varphi_1\xi_1=\eta_1$. Then 
$$ \gamma_1\varphi _1 \xi_1=\gamma_1 \eta_1\stackrel{\xi_1=\eta_1}{=}\gamma_1 \xi_1\neq \eta_1$$ by choice of $\gamma$. Moreover, we have $ \gamma_2\varphi _2 \xi_2 \neq \eta_2$, because
$ \gamma_2\varphi _2 \xi_2 = \eta_2=\gamma_2\xi_2$ implies that $\varphi_2=\gamma_2^{-1}\gamma_2\varphi_2$ is contained in the stabilizer of $\xi_2$ which is a contradiction to the choice of $\varphi$. Hence $\gamma\varphi$ is the desired element.

If $\xi=\eta$ we  choose $\gamma=(\gamma_1,\gamma_2)\in\Gamma$ \st $\gamma_2\xi_2\ne \eta_2$ and apply the first case. \qed\\[3mm]
Using (\ref{Fregaxial})  we immediately obtain the following
\begin{cor}\label{goodpos}
For any regular axial $h\in\Gamma$ and $\xi\in L_{\Gamma_1}\times L_{\Gamma_2}\subseteq\Frand$ there exists $\gamma\in\Gamma$ \st $\gamma\xi\in\horF(h^+)$.
\end{cor}

We now fix a regular axial isometry  $h=(h_1,h_2)\in \Gamma$  and a base point $\xo=(\xo_1,\xo_2)\in\XX$.
The following important theorem implies that $\Flim$ can be covered by
finitely many $\Gamma$-trans\-lates of an appropriate open set in $\Frand$. 
\begin{thr}\label{minimal}
The Furstenberg limit set is minimal, i.e. $\Flim$ is the smallest non-empty, $\Gamma$-invariant closed subset of $\Frand$.
\end{thr}
\prf\  We first show that every non-empty, $\Gamma$-invariant closed subset of $\Frand$ contains either $h^+$ or $h^-$. Replacing $h$ by its inverse if necessary, it then suffices to prove that $\Flim=\overline{\Gamma\cdot h^+}$.

Let $A\subseteq \Frand$ be a  non-empty, $\Gamma$-invariant closed set, and $\xi=(\xi_1,\xi_2)\in A$. If $\xi\in \{h^+, h^-\}$, there is nothing to prove, so assume that there exist indices 
$i,j\in\{1,2\}$ \st $\xi_i\ne h_i^+$ and $\xi_j\ne h_j^-$. If  $\xi_1\notin\{h_1^+,h_1^-\}$, then -- since $\xi_2$ is different from at least one of the points $h_2^+$,  $h_2^-$ -- 
we have $\xi\in\horF(h^+)$ or $\xi\in \horF(h^-)$.  So $\lim_{n\to\infty}h^{-n}\xi=h^-$ or $\lim_{n\to\infty}h^{n}\xi=h^+$ and we conclude that  $h^+$ or $h^-$ belongs to $A$. 
The case $\xi_2\notin\{h_2^+,h_2^-\}$  is analogous. It therefore remains to consider the possibilities $\xi=(h_1^-, h_2^+)$ or $\xi =(h_1^+,h_2^-)$. In both cases $\xi$ is contained in  
$L_{\Gamma_1}\times L_{\Gamma_2}$, so by Corollary~\ref{goodpos} there exists $\gamma\in\Gamma$ \st $\gamma\xi\in\horF(h^+)$. Then $\lim_{n\to\infty}h^{-n}\gamma\xi=h^-$ which proves that $h^-\in A$.

For the second part of the proof we are going to show the stronger statement that  $\Flim=\overline{\Gamma\cdot\xi}$ for any $\xi=(\xi_1,\xi_2)\in\Flim$. 

Let $\eta=(\eta_1,\eta_2)\in\Flim$ arbitrary. If $\eta=\xi$,  there is nothing to prove, if $\eta_1=\xi_1$ or $\eta_2=\xi_2$, then by Lemma~\ref{goodposition} there exists $\gamma=(\gamma_1,\gamma_2)\in\Gamma$ \st $\gamma_1\xi_1\neq\eta_1$ and $\gamma_2\xi_2\ne \eta_2$. Hence replacing $\xi$ by $\gamma\xi$ if necessary, we may assume that $\eta_1\ne\xi_1$ and $\eta_2\ne\xi_2$. Let $U_1\subset\rand_1$, $U_2\subset\rand_2$ be neighborhoods  of $\eta_1$, $\eta_2$   \st $\xi_1\notin U_1 $ and $\xi_2\notin U_2$, and choose $\tilde\eta\in\preim(\eta)\cap\Lim$. Then there exists a sequence $(\gamma_n)=\big((\gamma_{n,1},\gamma_{n,2})\big)\subset \Gamma$  \st $\gamma_n\xo\to \tilde\eta=(\eta_1,\eta_2,\theta)$ and $\gamma_n^{-1}\xo$ converges as $n\to\infty$. Since for $i=1,2$ $d_i(\xo_i,\gamma_{n,i}^{-1}\xo_i)=d_i(\gamma_{n,i}\xo_i,\xo_i)\to\infty$  we have $\gamma_{n,i}^{-1}\xo_i\to\zeta_i\in\rand_i$, $i=1,2$, and $\theta(\xo,\gamma_n^{-1}\xo)=\theta(\xo,\gamma_n\xo)\to \theta$ as $n\to\infty$. Hence $\lim_{n\to\infty}\gamma_n^{-1}\xo\to\tilde\zeta:=(\zeta_1,\zeta_2,\theta)\in\regrand$, and we put $\zeta:=\pi^F(\tilde\zeta)=(\zeta_1, \zeta_2)$. 

Moreover, we can assume $\zeta\in\horF(h^+)$, because otherwise, by Corollary~\ref{goodpos}, we find $\gamma\in\Gamma$ \st $\gamma\zeta\in\horF(h^+)$ and we can replace our sequence $(\gamma_n)$ by $(\gamma_n\gamma^{-1})$.

Let $T\gg1$, $\eps>0$ be arbitrary. Then Lemma~\ref{convergence} implies the existence of $N\in\NN$ \st for all $n\ge N$ and $t\in[0,T]$
$$d(\sigma_{\xo_i,\gamma_{n,i}\xo_i}(t),\sigma_{\xo_i,\gamma_{n,i} h_i^+}(t))=d(\sigma_{\gamma_{n,i}^{-1}\xo_i,\xo_i}(t),\sigma_{\gamma_{n,i}^{-1}\xo_i,h_i^+}(t))\le \frac{\eps}2$$
and $ d(\sigma_{\xo_i,\gamma_{n,i}\xo_i}(t),\sigma_{\xo_i,\eta_i}(t))\le\eps/2$.
Hence we conclude that as $n\to\infty$ $\gamma_{n,i}h_i^+\to\eta_i$ for $i=1,2$, in particular, there exists $\varphi=(\varphi_1,\varphi_2)\in\Gamma$ \st  $\varphi_i h_i^+\in U_i$ and $\varphi_i h_i\varphi_i^{-1}$ is rank one for $i=1,2$.

Assume first that $\xi\in\horF(\varphi h^-)$. Then there exists $N\in\NN$ \st for $n\ge N$  $((\varphi_1 h_1\varphi_1^{-1})^n\xi_1, (\varphi_2 h_2\varphi_2^{-1})^n\xi_2)=(\varphi h\varphi^{-1})^n\xi\in U_1\times U_2$.

If $\xi\notin\horF(\varphi h^-)$, Corollary~\ref{goodpos} implies the existence of $\gamma\in\Gamma$ \st $\gamma\xi\in\horF(\varphi h^-)$ and we conclude $(\varphi h\varphi^{-1})^n\gamma\xi\in U$ for $n$ sufficiently large.\qed 
\begin{thr}\label{Product}
 The regular geometric limit set $\Lim\cap\regrand$ is isomorphic to a product $\Flim\times P_\Gamma$, where $P_\Gamma\subseteq (0,\pi/2)$ denotes the set of  slopes of  regular limit points.
\end{thr}
\prf\  If $\tilde\xi\in \Lim\cap\regrand$, then $\pi^F(\tilde\xi)\in F_\Gamma$, and by definition of $P_\Gamma$ the slope of $\tilde\xi$ belongs to $P_\Gamma$. 

Conversely, let $\eta=(\eta_1,\eta_2)\in F_\Gamma$ and $\theta\in P_\Gamma$. We have to show that $\tilde\eta:=(\eta_1,\eta_2,\theta)\in\Lim$. By definition
of $P_\Gamma$ and Lemma~\ref{prodtopology} there
exists a sequence $(\gamma_n)\subset\Gamma$  \st $\theta_n:=\theta(\xo,\gamma_n\xo)$ converges to $\theta$ as $n\to\infty$. 
Moreover, by compactness of $\rand_1\times\rand_2$ a subsequence of $(\gamma_n\xo)$ converges to a point  $\tilde \xi\in\Lim\cap\regrand$ of slope $\theta$. Put $\xi:=\pi^F(\tilde\xi)$, and notice that $\tilde\eta\in\regrand$ is  the unique point in $(\pi^F)^{-1}(\eta)$ of slope $\theta$.

By {\rm Theorem}~\ref{minimal}  $\,F_\Gamma=\overline{\Gamma \cdot \xi}\,$ is a
minimal closed set under the action of $\Gamma$, hence 
$$\eta\in\overline{\Gamma\cdot\xi}=\pi^F(\overline{\Gamma\at\tilde\xi})\,.$$
Since the action of $\Gamma$ on the geometric boundary does not change
 the slope of a point, we conclude that  the closure of $\Gamma\cdot \tilde\xi$ contains $\tilde\eta$. In particular
$\tilde\eta \in\overline{\Gamma\at\tilde\xi}\subseteq\Lim.\ $\qed

\section{Density of  regular axial isometries}

In this section we will make the stronger assumption that $\Gamma\subset\is(\XX_1)\times\is(\XX_2)$ acts properly discontinuously on the product $\XX$ of two locally compact 
Hadamard spaces $\XX_1$, $\XX_2$ and contains two isometries $g=(g_1,g_2)$ and $h=(h_1,h_2)$ \st $g_1$ and $h_1$ are independent rank one elements of $\Gamma_1$ and $g_2$, $h_2$ are independent 
rank one elements in $\Gamma_2$. Recall that an isometry $h=(h_1,h_2)\in\is(\XX_1)\times\is(\XX_2)$ is called regular axial if $h_1$ and $h_2$ are rank one elements. Its attractive fixed point 
is denoted $\widetilde{h^+}\in\regrand$ and we put $h^+:=\pi^F(\widetilde{h^+})=(h_1^+, h_2^+)$. Moreover, for $h=(h_1,h_2)$ regular axial and $i\in\{1,2\}$ we denote $l_i(h_i)$ the translation length of $h_i$ in $\XX_i\,.\ $ The {\hl limit cone of $\Gamma$} is defined by
$$ \ell_\Gamma:=\overline{\{\arctan\big(l_2(g_2)/l_1(g_1)\big):g=(g_1,g_2)\in\Gamma\quad \mbox{regular axial}\}}\,.$$ 
We fix a base point $\xo=(\xo_1,\xo_2)\in\XX$. 
The following proposition is a key ingredient in the proofs.
\begin{prp}\label{gettingregularaxials}
Suppose $g=(g_1,g_2)$ and $h=(h_1,h_2)\in \Gamma$ are regular axial isometries \st $g_i$ and $h_i$ are independent in $\Gamma_i$ for $i=1,2$. Let $(\gamma_n)\subset\Gamma$ be a sequence 
\st  $\gamma_n\xo$ and $\gamma_n^{-1}\xo$ converge to points in $\regrand$ as $n\to\infty$. Then given arbitrarily small distinct neighborhoods $W_i(h^+), W_i(h^-)\subset\ganz_i$ of $h_i^+$, 
$h_i^-$, $i=1,2$, there exist $N\in\NN$, $\alpha\in\{h^N, h^Ng^N,h^Ng^{-N}\}$  and $\beta\in\{h^{-N},h^{-N}g^N, h^{-N}g^{-N}\}$ \st $\varphi_n:=\alpha\gamma_n\beta^{-1}$ satisfies 
$\varphi_{n}\xo\in W_1(h^+)\times W_2(h^+)$ and  $\varphi_{n}^{-1}\xo\in W_1(h^-)\times W_2(h^-)$ for $n$ sufficiently large.
\end{prp}
\prf\ For $i=1,2$ and $\eta\in\{g^-,g^{+},h^-,h^{+}\}$ let $W_i(\eta)\subset\ganz_i$ be an arbitrary, sufficiently small neighborhood of $\eta_i^+\in\rand_i$ with $\xo_i\notin W_i(\eta)$ \st all 
$W_i(\eta)$ are pairwise disjoint in $\ganz_i$. According to Lemma~\ref{dynrankone} (c) there exists  a constant   $N\in\NN$ \st for all $\gamma\in\{g,g^{-1},h,h^{-1}\}$ and $i\in\{1,2\}$
\begin{equation}\label{factorpingpong}  
\gamma_i^N\big(\ganz_i\setminus  W_i(\gamma^{-})\big)\subseteq W_i(\gamma^+)\,.
\end{equation}
Denote $F\subset\XX$ the finite set of points $\{\xo,h^{\pm N}\xo,  h^{\pm N}g^{\pm N}\xo\}$ in $\XX$. 
Since $\gamma_{n,i}\xo_i$ converges to a point $\xi_i\in\rand_i$, $i=1,2$, given arbitrary neighborhoods $U_1\subset\ganz_1$, $U_2\subset\ganz_2$ of $\xi_1$, $\xi_2$, there exists $N_+\in\NN$ \st for all $n>N_+$ and every $x\in F$ we have $\gamma_{n} x\in U_1\times U_2$.  
Using the fact that  $\ganz_i=\big(\ganz_i\setminus W_i(g^-)\big) \cup \big(\ganz_i\setminus W_i(g^+)\big)$, and choosing the neighborhoods $U_i$ of $\xi_i$, $i\in\{1,2\}$,  sufficiently small, we may assume that one of the following six  possibilities occurs for all $n>N_+$ and every $x=(x_1,x_2)\in F$:
\begin{enumerate}
\item[1.]Case: $\gamma_{n,1}x_1\in \ganz_1\setminus W_1(h^-)$ and $\gamma_{n,2}x_2\in\ganz_2\setminus W_2(h^-)$\\
Then by~(\ref{factorpingpong}) $\,h^{N} \gamma_nx\in W_1(h^+)\times W_2(h^+)$.
\item[2.]Case: $\gamma_{n,1}x_1\in W_1(h^-)$ and $\gamma_{n,2}x_2\in W_2(h^-)$\\
Since $W_i(h^-)\subset \ganz_i\setminus W_i(g^-)$, $i=1,2$, we have again by ~(\ref{factorpingpong}) $\,g^N\gamma_nx\in W_1(g^+)\times W_2(g^+)$. Hence we are in case 1 for $g^N \gamma_nx$, so $h^N g^N\gamma_nx\in W_1(h^+)\times W_2(h^+)$.
\item[3.]Case:  $\gamma_{n,1}x_1\in W_1(h^-)$ and $\gamma_{n,2}x_2\in\ganz_2\setminus \big(W_2(h^-)\cup W_2(g^-)\big)$\\
Then $g^N\gamma_nx\in W_1(g^+)\times W_2(g^+)$, which yields $h^Ng^N\gamma_nx \in W_1(h^+)\times W_2(h^+)$.
\item[4.]Case:  $\gamma_{n,1}x_1\in W_1(h^-)$ and $\gamma_{n,2}x_2\in\ganz_2\setminus \big(W_2(h^-)\cup W_2(g^+)\big)$\\
Then $g^{-N}\gamma_nx\in W_1(g^-)\times W_2(g^-)$, which gives $h^Ng^{-N}\gamma_nx \in W_1(h^+)\times W_2(h^+)$.
\item[5.]Case: $\gamma_{n,1}x_1\in \ganz_1\setminus \big(W_1(h^-)\cup W_1(g^-)\big)$ and $\gamma_{n,2}x_2\in W_2(h^-)$\\
Similarly to case 3 we obtain $h^N g^N\gamma_nx\in W_1(h^+)\times W_2(h^+)$.
\item[6.]Case:  $\gamma_{n,1}x_1\in \ganz_1\setminus  \big(W_1(h^-)\cup W_1(g^+)\big)$ and $\gamma_{n,2}x_2\in W_2(h^-)$\\
As in case 4 we get $h^Ng^{-N}\gamma_nx\in W_1(h^+)\times W_2(h^+)$.
\end{enumerate}
So we have shown the existence of  $\alpha\in\{ h^N, h^Ng^N, h^Ng^{-N}\}$ \st for all $n>N_+$ and every  $x\in F$  $\alpha \gamma_n x \in W_1(h^+)\times W_2(h^+)$.

With a similar case by case treatment  we get   $N_-\in\NN$ and $\beta\in \{h^{-N}, h^{-N} g^N, h^{-N}g^{-N}\}$ \st $\beta\gamma_n^{-1}x\in W_1(h^-)\times W_2(h^-)$ for all $n>N_2$ and all $x\in F$. Then, putting $\varphi_n:=\alpha \gamma_n\beta^{-1}$,  the claim holds for all $n>\max\{N_+,N_-\}$.\qed\\[3mm]
The following theorem relates the limit cone to $P_\Gamma$.
\begin{thr}\label{propertiesofPGamma}
If $\Gamma$ contains two regular axial isometries which project to independent rank one elements in each factor  then $P_\Gamma=\ell_\Gamma\cap (0,\pi/2)$.  Moreover, $\ell_\Gamma$ is either a point or an interval.
\end{thr}
\prf\ We first prove $\ell_\Gamma\cap (0,\pi/2)\subseteq P_\Gamma$: If $g_n=(g_{n,1},g_{n,2})$ is a sequence of regular axial isometries \st $\arctan\big(l_2(g_{n,2})/l_1(g_{n,1})\big)$ converges to $\theta\in (0,\pi/2)$, we choose 
$$k_n\ge 2n \max\{d_i(\xo_i,\Ax(g_{n,i}))/l_i(g_{n,i}):i=1,2\}$$ and put
$\gamma_n:=g_n^{k_n}$. From
$$ k_n l_i(g_{n,i})\le d_i(\xo_i,\gamma_{n,i}\xo_i)\le 2 d_i(\xo_i,\Ax(g_{n,i}))+k_nl_i(g_{n,i}) \le k_n l_i(g_{n,i})(1+1/n)$$
we  get 
\be
\tan \theta & = &\lim_{n\to\infty}\left(\frac{l_2(g_{n,2})}{l_1(g_{n,1})}\cdot (1+\frac1{n})\right) \ge \lim_{n\to\infty}\frac{d_2(\xo_2,\gamma_{n,2}\xo_2)}{d_1(\xo_1,\gamma_{n,1}\xo_1)}\,,\\
\tan \theta & = &\lim_{n\to\infty}\left(\frac{l_2(g_{n,2})}{l_1(g_{n,1})}\cdot \frac{n}{n+1}\right) \le \lim_{n\to\infty}\frac{d_2(\xo_2,\gamma_{n,2}\xo_2)}{d_1(\xo_1,\gamma_{n,1}\xo_1)}\,,
\ee
hence the claim. 

Let's prove the inclusion $P_\Gamma\subseteq \ell_\Gamma\cap (0,\pi/2)$. Denote $g=(g_1, g_2)$, $h=(h_1,h_2)\in\Gamma$ two regular axial isometries as in Proposition~\ref{gettingregularaxials}. For $\eta\in \{g^-,g^{+},h^-,h^{+}\}$ and  $i\in\{1,2\}$  let $U_i(\eta)\subset \ganz_i$ be a small neighborhood of $\eta_i$ with $\xo_i\notin U_i(\eta)$ \st all $U_i(\eta)$ are pairwise disjoint. Upon taking smaller neighborhoods, Lemma~\ref{joinrankone} provides a constant  $c>0$  \st for $i\in\{1,2\}$ any pair of points  in distinct neighborhoods can be joined by a rank one geodesic $\sigma_i\subset\XX_i$ with $d(\xo_i,\sigma_i)\le c$.  Moreover, according to Lemma~\ref{elementsarerankone} for $i\in\{1,2\}$ and $\eta\in \{h^-,h^{+}\}$ there exist neighborhoods $W_i(\eta)\subseteq U_i(\eta)$ of $\eta_i$  \st every $\gamma=(\gamma_1,\gamma_2)\in\Gamma$ with $\gamma_i\xo_i\in W_i(h^+)$ and  $\gamma_i^{-1}\xo\in W_i(h^-)$, $i=1,2$,   is regular axial with $\gamma_i^+\in U_i(h^+)$ and $\gamma_i^-\in U_i(h^-)$, $i=1,2$. 

Now let $\theta\in P_\Gamma$. By definition there exists a sequence $(\gamma_n)=\big((\gamma_{n,1},\gamma_{n_2})\big)\subset\Gamma$ \st $d_2(\xo_2,\gamma_{n,2}\xo_2)/d_1(\xo_1,\gamma_{n,1}\xo_1)\to\tan \theta$, $\gamma_{n,1}\xo_1\to \xi_1$, $\gamma_{n,2}\xo_2\to\xi_2$ as $n\to\infty$. Passing to a subsequence if necessary, we can assume that $\gamma_{n}^{-1}\xo\to \tilde\zeta=(\zeta_1,\zeta_2,\theta)\in\regrand$  as $n\to\infty$.  By Proposition~\ref{gettingregularaxials} there exist $N_0\in\NN$,  a finite set $\Lambda\subset\Gamma$ and $\alpha$, $\beta\in\Lambda$ \st for all $n>N_0$
$$ \alpha\gamma_n\beta^{-1}\xo\in W_1(h^+)\times W_2(h^+)\ \mbox{and}\quad \beta\gamma_n^{-1}\alpha^{-1}\xo\in W_1(h^-)\times W_2(h^-)\,.$$
Put $\varphi_n:=\alpha\gamma_n\beta^{-1}$, $n\in\NN$, and  
$L:=\max\{d_i(\xo_i,\lambda_i\xo_i): i\in\{1,2\}, \lambda=(\lambda_1,\lambda_2)\in\Lambda \}$. Using  the triangle inequality we estimate  for $i=1,2$
\begin{equation}\label{samedistance}
|d_i(\xo_i,\varphi_{n,i}\xo_i)-d_i(\xo_i,\gamma_{n,i}\xo_i)|\le 2L\,.
\end{equation}
Moreover, by choice of the sets $W_i(h^\pm)\subseteq U_i(h^\pm)$ we know that for $n>N_0$  $\varphi_n$ is regular axial with $\varphi_n^+\in U_1(h^+)\times U_2(h^+)$ and $\varphi_n^-\in U_1(h^-)\times U_2(h^-)$. So Lemma~\ref{joinrankone} shows that for $n>N_0$ there exists 
 $x_{n,i}\in \Ax(\varphi_{n,i})$ \st $d_i(\xo_i, x_{n,i})\le c$, $i=1,2$. We conclude 
$$ l_i(\varphi_{n,i})\le d_i(\xo_i,\varphi_{n,i}\xo_i)\le l_i(\varphi_{n,i})+2c\,,\ i=1,2\,
,$$
which -- together with (\ref{samedistance}) --  implies that $\tan\theta=\lim_{n\to\infty} l_2(\varphi_{n,2})/l_1(\varphi_{n,1})$. 
 
Let's prove the last assertion following the lines of the proof of Proposition 2.4 in \cite{DalboKim}: If $\ell_\Gamma$ is a point,  there is nothing to prove.  Otherwise we will show that for $\theta,\theta'\in \{\arctan\big(l_2(g_2)/l_1(g_1)\big):g=(g_1,g_2)\in\Gamma\ \,\mbox{regular axial}\}$, $\theta <\theta'$, we have $[\theta,\theta']\subseteq \ell_\Gamma$.
Fix $\gamma=(\gamma_1,\gamma_2)$, $\varphi=(\varphi_1,\varphi_2)$ regular axial \st $\tan \theta=l_2(\gamma_2)/l_1(\gamma_1)$ and $\tan \theta'=l_2(\varphi_2)/l_1(\varphi_1)$.

Recall that $g$, $h\in\Gamma$ are two regular axial isometries projecting to independent rank one elements. If $\gamma_1$, $\varphi_1$ are not independent, then by Lemma~\ref{dynrankone} (c) there exist $N\in\NN$, $\alpha, \beta \in\{h^N, h^{-N}, g^N, g^{-N}\}$ with $\alpha\ne \beta$ \st $\alpha_1 \gamma_1\alpha_1^{-1}$ and $\beta_1\varphi_1 \beta_1^{-1}$ are independent. Using the fact that the translation length is invariant by conjugation and upon replacing $\gamma$ by $\alpha \gamma\alpha^{-1}$ and $\varphi$ by $\beta \varphi\beta^{-1}$ if necessary, we may assume that $\gamma_1$ and $\varphi_1$ are independent rank one elements of $\Gamma_1$. 

Now either $\gamma_2$ and $\varphi_2$ are independent, or, after replacing $\gamma$, $\varphi$ by its inverse if necessary, we have $\gamma_2^+=\varphi_2^+$. By Proposition~\ref{combgeomlength} and Lemma~\ref{combgeomlen} there exist $N\in\NN$ and $C>0$ \st for $i\in\{1,2\}$ and all $n,m\in\NN\setminus\{0\}$  
$$ \big| l_i(\gamma_i^{Nn}\varphi_i^{Nm})-Nn\;l_i(\gamma_1)-Nm\;l_i(\varphi_1)\big|\le C\,.$$
Hence
$$\lim_{k\to\infty}\frac{l_2(\gamma_2^{Nnk}\varphi_2^{Nmk})}{l_1(\gamma_1^{Nnk}\varphi_1^{Nmk})} =\frac{n\;l_2(\gamma_2)+m\;l_2(\varphi_2)}{n\;l_1(\gamma_1)+m\;l_1(\varphi_1)}\,,$$
so we have 
$$ \arctan\left(\frac{l_2(\gamma_2)+q\;l_2(\varphi_2)}{l_1(\gamma_1)+q\;l_1(\varphi_1)}\right)\in\ell_\Gamma$$
for every positive rational number $q\in\QQ$ and we conclude $[\theta,\theta']\in\ell_\Gamma$.\qed\\[3mm]
In order to prove Theorem~D from the introduction, we will need an important definition as a substitute for  the more familiar notion of $\Gamma$-duality used e.g. in \cite{MR1383216} and \cite{GR/0809.0470} when dealing with only one factor.
\begin{df}
Two points $\xi=(\xi_1,\xi_2)$, $\eta=(\eta_1,\eta_2)\in \Frand$ are called {\hd $\Gamma$-related} if for any neighborhoods $U_1$, $V_1\subset\ganz_1$ of $\xi_1$, $\eta_1$ and all neighborhoods $U_2$, $V_2\subset\ganz_2$ of $\xi_2$, $\eta_2$ there exists $\gamma=(\gamma_1,\gamma_2)\in\Gamma$ 
\st  for $i\in\{1,2\}$
$$ \gamma_i(\ganz_i\setminus U_i)\subset V_i\,,\quad \gamma_i^{-1}(\ganz_i\setminus V_i)\subset U_i\,.$$
We will denote $\rel(\xi)$ the set of points in $\Frand$ which are $\Gamma$-related to $\xi$.
\end{df}
Notice that for any $\xi\in\Frand$ the set $\rel(\xi)$ is closed with respect to the topology of $\Frand$. Moreover, if $\eta\in\rel(\xi)$, then $\eta_1$ is $\Gamma_1$-dual to $\xi_1$ and $\eta_2$ is $\Gamma_2$-dual to $\xi_2$. % The converse clearly does not hold in general.

The importance of the notion lies in the following. If $\tilde h^+$, $\tilde h^-$ denote the attractive and repulsive fixed point of a regular axial isometry $h=(h_1,h_2)\in\Gamma$, then $h^+= \pi^F(\tilde h^+)$ and $h^-= \pi^F(\tilde h^-)$ are $\Gamma$-related by  Lemma~\ref{dynrankone} (c). Conversely, if $\xi=(\xi_1,\xi_2)$, $\eta=(\eta_1,\eta_2)\in\Frand$ are $\Gamma$-related, then by definition there exists a sequence $(\gamma_n)=\big((\gamma_{n,1},\gamma_{n,2})\big)\subset\Gamma$ \st for $i\in\{1,2\}$ we have $\gamma_{n,i}\xo_i\to\eta_i$ and $\gamma_{n,i}^{-1}\xo_i\to\xi_i$ as $n\to\infty$. Hence if $\xi_i$ can be joined to $\eta_i$ by a rank one geodesic for $i=1,2$, then in view of  Lemma~\ref{elementsarerankone} $\gamma_n$ is regular axial for $n$ sufficiently large and satisfies 
$$\gamma_{n,i}^+\to\eta_i\quad\mbox{and}\qquad \gamma_{n,i}^{-}\to\xi_i\quad\mbox{as}\  n\to\infty$$
for $i\in\{1,2\}$ . Denote $\Delta \subset\Frand\times\Frand$ the set 
$$\Delta:=\{(\xi,\eta)\in \Frand\times\Frand : \xi_1= \eta_1\ \mbox{or}\  \ \xi_2=\eta_2\}\,.$$

Using the above definition, we are now able to prove the following statement which is Theorem~D from the introduction and can be viewed as a  strong topological version of the double ergodicity property of Poisson boundaries due to Burger-Monod (\cite{MR1911660}) and Kaimanovich (\cite{MR2006560}). 
\begin{thr}\label{densityofaxials}
If $\Gamma$ contains two regular axial isometries projecting to independent rank one elements in each factor then the set of pairs of fixed points $(\gamma^+,\gamma^-)\subset\Frand\times \Frand$ of regular axial isometries  $\gamma\in\Gamma$ is dense in $\big(F_\Gamma
\times F_\Gamma\big)\setminus \Delta$.
\end{thr}
\prf\ Denote $g=(g_1,g_2)$ and  $h=( h_1, h_2)\in\Gamma$ two regular axial isometries \st for $i\in\{1,2\}$ $g_i$ and $h_i$ are independent. In view of the paragraph preceding the theorem we 
first prove that any two distinct points  in $\{g^-,g^+,h^-, h^+\}$ are $\Gamma$-related. 

For $\eta\in \{g^-,g^{+},h^-,h^{+}\}$ and $i\in\{1,2\}$ let $U_i(\eta)\subset \ganz_i$ be an arbitrary, sufficiently small neighborhood of $\eta_i$ with $\xo_i\notin U_i(\eta)$ \st all $U_i(\eta)$ are pairwise disjoint. According to Lemma~\ref{dynrankone} (c) there exists  a constant   $N\in\NN$ \st for all $\gamma\in\{g,g^{-1},h,h^{-1}\}$ and $i\in\{1,2\}$
\begin{equation}\label{factorpingpongg}  
\gamma_i^N\big(\ganz_i\setminus U_i(\gamma^{-})\big)\subseteq U_i(\gamma^+)\,.
\end{equation}
Let $\gamma,\varphi \in\{g,g^{-1},h,h^{-1}\}$, $\varphi\ne \gamma$.  Using the fact that either $\varphi=\gamma^{-1}$ or $\gamma_i$, $\varphi_i$ are independent for $i=1,2$  (\ref{factorpingpongg}) implies 
\be
&&\gamma_i^N \underbrace{\varphi_i^{-N}\big(\ganz_i\setminus U_i(\varphi^+)\big)}_{\subset U_i(\varphi^{-})\subset \ganz_i\setminus U_i(\gamma^{-})}\subset U_i(\gamma^+)\quad \mbox{and}\\
&&(\gamma_i^N\varphi_i^{-N})^{-1}\big(\ganz_i\setminus U_i(\gamma^+)\big)\subset\varphi_i^N \big(U_i(\gamma^{-})\big)\subset U_i(\varphi^+)\,
\ee
 for $i\in\{1,2\}$. Hence $\varphi^+\in\rel(\gamma^+)$.
 
 Next we will show that any $\xi=(\xi_1,\xi_2)\in\Flim$ with $\xi_i\notin\{g_i^-,g_i^+,h_i^-,h_i^+\}$, $i=1,2$, is $\Gamma$-related to any point in $\{g^-,g^+,h^-,h^+\}$.  
For $\zeta\in \{\xi, g^-,g^{+},h^-,h^{+}\}$ and $i\in\{1,2\}$ let $U_i(\zeta)\subset \ganz_i$ be an arbitrary, sufficiently small neighborhood of $\zeta_i$ with $\xo_i\notin U_i(\zeta)$  \st all $U_i(\zeta)$ are pairwise disjoint.  By  Lemma~\ref{elementsarerankone} there exist neighborhoods $W_i(\zeta)\subseteq U_i(\zeta)$, $\zeta\in \{\xi, g^-,g^{+},h^-,h^{+}\}$, \st every $\gamma_i\in\Gamma_i$ with $\gamma_i\xo_i\in W_i(\zeta)$, $\gamma_i^{-1}\xo_i\in W_i(\eta)$, $\eta\in\{\xi, g^-,g^{+},h^-,h^{+}\}\setminus\{\zeta\}$,  is rank one with $\gamma_i^+\in U_i(\zeta)$ and $\gamma_i^-\in U_i(\eta)$. 

Since $\xi\in\Flim$, there exists a sequence $(\gamma_n)=\big((\gamma_{n,1},\gamma_{n,2})\big)\subset\Gamma$ \st $\gamma_{n,1}\xo_1\to \xi_1$, $\gamma_{n,2}\xo_2\to\xi_2$. Upon passing to  a subsequence if necessary we may assume that $\gamma_{n,1}^{-1}\xo_1\to \zeta_1\in\rand_1$ and $\gamma_{n,2}^{-1}\xo_2\to\zeta_2\in\rand_2$. By Proposition~\ref{gettingregularaxials} there exist $N, N_0\in\NN$ and $\beta\in\{h^{-N}, h^{-N}g^N, h^{-N} g^{-N} \}$ \st for all $n>N_0$  $\gamma_n\beta^{-1}\xo\in W_1(\xi)\times W_2(\xi)$ and $\beta\gamma_n^{-1}\xo\in W_1(h^-)\times W_2(h^-)$. 
By  Lemma~\ref{elementsarerankone} we conclude that for $n>N_0$ the isometry $\gamma_n\beta^{-1}$ is regular axial with $(\gamma_n\beta^{-1})^+\in U_1(\xi)\times U_2(\xi)$ and $(\gamma_n\beta^{-1})^-\in  U_1(h^-)\times U_2(h^-)$. This implies that $\xi\in\rel(h^-)$ and by symmetry 
\begin{equation}\label{relatedtofourpoints}
\xi\in\rel(g^-)\cap \rel(g^+)\cap\rel(h^-)\cap\rel(h^+)\,.
\end{equation}

Next we let $\xi=(\xi_1,\xi_2),\eta=(\eta_1,\eta_2)\in\Flim$  \st for $i\in\{1,2\}$  we have $\xi_i, \eta_i\notin\{g_i^-,g_i^+,h_i^-,h_i^+\}$ and $\xi_i\ne\eta_i$. As above, for $\zeta\in \{\xi,\eta,  h^-\}$ and $i\in\{1,2\}$ let $U_i(\zeta)\subset \ganz_i$ be an arbitrary, sufficiently small neighborhood of $\zeta_i$ with $\xo_i\notin U_i(\zeta)$ \st all $U_i(\zeta)$ are pairwise disjoint.  
By the arguments in the previous paragraph there exists a regular axial isometry $\varphi\in\Gamma$ with $\varphi^+\in U_1(\xi)\times U_2(\xi)$ and $\varphi^-\in U_1(h^{-})\times U_2(h^{-})$. In particular, $\varphi_i$ and $g_i$ are independent for $i=1,2$. Replacing $h$ by $\varphi$ in (\ref{relatedtofourpoints}) we know that $\eta\in \rel(g^-)\cap \rel(g^+)\cap\rel(\varphi^-)\cap\rel(\varphi^+)$, in particular $\eta\in\rel(\varphi^+)$. So using the fact that $\eta_i$ can be joined to $\varphi_i^+$ by a rank one geodesic in $\XX_i$ for $i=1,2$, given small neighborhoods $U_i(\varphi^+)\subseteq U_i(\xi)$ for $i\in\{1,2\}$, there exists $\gamma\in\Gamma$ regular axial  with $\gamma^+\in U_1(\varphi^+)\times U_2(\varphi^+)\subseteq U_1(\xi)\times U_2(\xi)$ and $\gamma^-\in U_1(\eta)\times U_2(\eta)$. \qed

\section{The exponent of growth for a given slope}\label{ExpGrowthSlope}

For the remainder of the article  $\XX$ is a product of locally compact 
Hadamard spaces $\XX_1$, $\XX_2$, $\xo=(\xo_1, \xo_2)$ a fixed base point, and $\Gamma\subset\is(\XX_1)\times\is(\XX_2)$ a discrete group which contains two isometries $g=(g_1,g_2)$ and $h=(h_1,h_2)$ \st $g_i$ and $h_i$ are independent rank one elements of $\Gamma_i$ for $i=1,2$.  
In this section we want to describe the map which assigns to each slope $\theta\in [0,\pi/2]$ the exponential growth rate of orbit points of $\Gamma$ in $\XX$ with a prescribed slope $\theta$.  Recall the notation introduced in Section~\ref{prodHadspaces} and put for $x,y\in\XX$, $\theta\in[0,\pi/2]$, $\eps>0$ 
$$ \Gamma(x,y;\theta,\eps):=\{\gamma\in\Gamma: \gamma y\ne x\quad \mbox{and}\ \ |\theta(x,\gamma y)-\theta | <\eps \}\,.$$
For the definition of the exponential growth rate we introduce the following partial sum of the Poincar{\'e} series for $\Gamma$: $\;$  For $s>0$ we put
$$Q^{s,\eps}_{\theta}(x,y)=\sum_{\gamma\in
  \Gamma(x,y;\theta,\eps)}
  e^{-s d(x,\gamma y)}$$
and denote $\delta_\theta^\eps(x,y)$ its {\hl critical exponent}, i.e. the   unique
  real number \st $Q^{s,\eps}_{\theta}(x,y)$ converges if
  $s>\delta_\theta^\eps(x,y)$ and diverges if
  $s<\delta_\theta^\eps(x,y)$. It is clear that for any $\eps>0$ we have $\delta_\theta^\eps(x,y)\le \delta(\Gamma)$, the critical exponent of the Poincar{\'e} series. Unfortunately, unlike in 
the case of $\delta(\Gamma)$, where the summation is over all elements in $\Gamma$, this number may depend on $x$ and $y$. If $\eps> \pi/2$ then the summation above is over all $\gamma\in \Gamma$ 
with $\gamma y \ne  x$. By discreteness of $\Gamma$ we have $\gamma y = x$ for only finitely many $\gamma\in\Gamma$, hence for $\eps>\pi/2$ we have $\delta_\theta^\eps(x,y)= \delta(\Gamma)$. 

For $n\in\NN$ we define 
$$N_\theta^\eps(x,y;n):=\#\{ \gamma\in\Gamma\;:\, n-1< d(x,\gamma y)\le n\,,\ |\theta(x,\gamma y)-\theta|<\eps\}\,,$$ 
which can be interpreted as an orbit counting function for orbit points of slope $\eps$-close to $\theta$.  Although the proof of the following lemma is standard, we include it here for the convenience of the reader.
\begin{lem}\label{defbylimsup}
We have 
$$\delta_{\theta}^{\eps}(x,y)=\limsup_{n\to\infty}\frac{\log
  N_{\theta}^{\eps}(x,y;n)}{n}\,.$$
 \end{lem}
\prf\  We clearly have
$$ Q_\theta^{s,\eps}(x,y)=\sum_{n=1}^\infty \sum_{\begin{smallmatrix}{\gamma\in\Gamma(x,y;\theta,\eps)}\\{n-1<d(x,\gamma y)\le n}\end{smallmatrix}} e^{-s d(x,\gamma y)}\,,$$
hence 
$$\sum_{n=1}^\infty e^{-sn} N_\theta^\eps(x,y;n)\le Q_\theta^{s,\eps}(x,y)\le \sum_{n=1}^\infty e^{-s(n-1)} N_\theta^\eps(x,y;n)=e^s \sum_{n=1}^\infty e^{-sn} N_\theta^\eps(x,y;n)\,.$$
Moreover, we can write 
$$ e^{-sn} N_\theta^\eps(x,y;n)= \left( e^{-s+\frac{\log N_\theta^\eps (x,y;n)}{n}}\right)^n\,,$$
so  finding an estimate for the term in the bracket independent of $n$ will allow us to compare $Q_\theta^{s,\eps}(x,y)$ to a geometric series.

Suppose first that $s> \limsup_{n\to\infty}\frac{\log
  N_{\theta}^{\eps}(x,y;n)}{n}$. Then there exists $N\in\NN$ \st for any $n\ge N$
  $$ \frac{\log N_{\theta}^{\eps}(x,y;n)}{n}<s\,$$
  and we estimate
  $$ Q_\theta^{s,\eps}(x,y)\le e^s \left( \sum_{n=1}^{N-1} e^{-sn} N_\theta^\eps(x,y;n)+\sum_{n=N}^\infty  \left( e^{-s+\frac{\log N_\theta^\eps (x,y;n)}{n}}\right)^n \right)\,.$$ 
  The first sum is finite, and the second term converges because the number inside the brackets is strictly smaller than $1$ for all $n\ge N$.
  
  If $s< \limsup_{n\to\infty}\frac{\log N_{\theta}^{\eps}(x,y;n)}{n}$, there exists a strictly increasing sequence $(n_k)\subset \NN$ \st $\lim_{k\to\infty} \frac{ \log N_\theta^\eps (x,y;n_k)}{n_k}>s$. 
In particular  there exists $N\in\NN$ \st $\frac1{n_k} \log N_\theta^\eps (x,y;n_k)>s$ for any $k\ge N$. Moreover, since $n_k\ge k$ for all $k$, we have 
  $$ Q_\theta^{s,\eps}(x,y)\ge \sum_{k=N}^\infty  \left( e^{-s+\frac{\log N_\theta^\eps (x,y;n_k)}{n_k}}\right)^{n_k}\ge  \sum_{k=N}^\infty  \left( e^{-s+\frac{\log N_\theta^\eps (x,y;n_k)}{n_k}}\right)^{k}\,,$$
which shows that  $Q_\theta^{s,\eps}(x,y)$ diverges.\qed
\begin{df}
The number $\delta_{\theta}(\Gamma):=\liminf_{\eps\to
  0}\delta_{\theta}^{\eps}(o,o)$ is called the {\hd exponent of
  growth of $\Gamma$ of slope $\theta$}.
\end{df}
The following lemma shows that this number  $\delta_\theta(\Gamma)$  does not depend on the choice of arguments of $\delta_{\theta}^\eps$.
\begin{lem}\label{indearg}
For $x,y\in\XX$ arbitrary we have $\ \liminf_{\eps\to
  0}\delta_{\theta}^{\eps}(x,y)=\delta_{\theta}(\Gamma)$.
\end{lem}
\prf\  Fix $\theta\in [0,\pi/2]$ and set 
$$H_\theta:=\left(\begin{array}{c}\cos\theta\\\sin\theta\end{array}\right)\,.$$ 
We first note that for any $x=(x_1,x_2)$ ,$y=(y_1,y_2)\in \XX$ and $\gamma=(\gamma_1,\gamma_2)\in\Gamma$  the equality
\begin{equation}\label{scalprodcos}
\langle H(x,\gamma y),H_\theta\rangle= d(x,\gamma y)\cdot \cos\big(\theta(x,\gamma y)-\theta\big)
\end{equation}
holds. Using
$$H(x,\gamma y)-H(\xo,\gamma\xo)= \left(\begin{array}{c} d_1(x_1,\gamma_1 y_1)-d_1(\xo_1,\gamma_1\xo_1)\\ d_2(x_2,\gamma_2 y_2)-d_2(\xo_2,\gamma_2\xo_2)\end{array}\right) \,,$$
 setting $c:=d(x,\xo)+d(y,\xo)$ and recalling that both $\sin\theta$ and $\cos\theta$ belong to the interval $[0,1]$ we further have
  $$| \langle H(x,\gamma y)-H(\xo,\gamma\xo),H_\theta\rangle |\le 4c\,. $$
In particular, we conclude
\be
\langle H(x,\gamma y),H_\theta\rangle &\ge & \langle H(\xo,\gamma\xo),H_\theta\rangle -4c = d(\xo,\gamma\xo)\cdot  \cos\big(\theta(\xo,\gamma\xo)-\theta\big) -4c\\
&\ge & d(x,\gamma y) \cdot  \cos\big(\theta(\xo,\gamma\xo)-\theta\big) -6c\,,
\ee
hence $$\cos\big(\theta(x,\gamma y)-\theta\big)\ge \cos\big(\theta(\xo,\gamma\xo)-\theta\big) -\frac{6c}{d(x,\gamma y)}\,.$$
This shows that given $\eps >0$, there exists $R\gg 1$ \st $d(x,\gamma y)>R$ and $|\theta(\xo,\gamma \xo)-\theta|<\frac\eps{2}$ implies
$ |\theta(x,\gamma y)-\theta|<\eps$.
A symmetric argument -- with the roles  of $(x,\gamma y)$ and $(\xo,\gamma\xo)$ exchanged -- ensures the existence of $R'\gg 1$ \st $d(x,\gamma y)>R'$ and $ |\theta(x,\gamma y)-\theta|<\eps$ implies $|\theta(\xo,\gamma \xo)-\theta|<2 \eps$. Summarizing, we know that for any $\eps >0$  there exists $R\gg 1$ \st for any $\gamma\in\Gamma$ with $d(x,\gamma y)>R$ we have the implications
$$ \gamma\in   \Gamma(\xo,\xo;\theta,\eps/2) \quad  \Longrightarrow\qquad   \gamma\in\Gamma(x,y;\theta,\eps) \quad \Longrightarrow\qquad  \gamma\in  \Gamma(x,y;\theta,2\eps).$$
Since by discreteness of $\Gamma$ there are only finitely many $\gamma\in\Gamma$ with $d(x,\gamma y)\le R$, we conclude that for any $\eps>0$ 
$$ \delta_\theta^{\eps/2}(\xo,\xo)\le \delta_\theta^\eps(x,y)\le \delta_\theta^{2\eps}(\xo,\xo).$$
Taking the limit inferior as $\eps$ tends to zero finishes the proof. \qed\\[3mm]
Notice that in the definition of $\delta_\theta(\Gamma)$ for $\theta\in (0,\pi/2)$ one may substitute 
$$\#\{ \gamma\in\Gamma\;:\, \gamma y\ne x\,,\ d(x,\gamma y)\le n\,,\ \Big| \frac{d_2(p_2(\gamma y),p_2(x))}{d_1(p_1(\gamma y),p_1(x))}-\tan\theta\Big|<\eps\}\,$$
in~(\ref{defbylimsup}) instead of $N_\theta^\eps(x,y;n)$. 
Furthermore, the following property holds:
\begin{lem}\label{nichtleer}
If $\Lim\cap \rand_\theta \ne \emptyset$, then $\delta_{\theta}(\Gamma)\ge 0$.
\end{lem}
\prf\  Suppose $\Lim\cap \rand_\theta\ne \emptyset$.
Then by Lemma~\ref{prodtopology} for any $\eps>0$ there
exist infinitely many $\gamma\in\Gamma$ \st $|\theta(\xo,\gamma\xo)-\theta|<\eps$. In particular
$$ \sum_{\gamma\in
  \Gamma(o,o;\theta,\eps)} 1 =
  Q^{\;0,\eps}_{\theta}(o,o) \quad\mbox{diverges}\,,$$
hence $\delta_{\theta}^\eps(o,o)\ge 0$. We conclude $\delta_{\theta}(\Gamma)=\liminf_{\eps\to 0}\delta_{\theta}^{\eps}(o,o)\ge 0$.\qed\\[3mm]
The following proposition states that the map $\theta\mapsto\delta_\theta(\Gamma)$ is upper semi-continuous.
\begin{prp}\label{upsemcon}
Let $(\theta_j)\subset [0,\pi/2]$ be a sequence converging to $\theta\in [0,\pi/2]$. Then
$$ \limsup_{j\to\infty}\,\delta_{\theta_j}(\Gamma)\le \delta_{\theta}(\Gamma)\,.$$
\end{prp}
\prf\  Let $\eps_0\in (0,\pi/2)$. Then $\theta_j\to\theta$ implies
$|\theta_j-\theta|<\eps_0/2$ for $j$ sufficiently large. Let
$\eps\in (0,\eps_0/2)$ and $\gamma\in\Gamma(\xo,\xo;\theta_j,\eps)$.  Then
$$|\theta(\xo,\gamma\xo)-\theta|<\eps+\eps_0/2<\eps_0\,,$$
hence for $j$ sufficiently large  $\Gamma(\xo,\xo:\theta_j,\eps)\subseteq \Gamma(\xo,\xo:\theta,\eps_0)$.
This shows 
$\delta_{\theta_j}^{\eps}(o,o)\le\delta_{\theta}^
{\eps_0}(o,o)$, and therefore
$\ \delta_{\theta_j}(\Gamma)=\liminf_{\eps\to
  0}\,\delta_{\theta_j}^{\eps}(o,o)\le\delta_{\theta}^{\eps_0}(o,o)\,.$\\[2mm]
We conclude \vspace{-2mm}
\begin{eqnarray*}
 \limsup_{j\to\infty} \delta_{\theta_j}(\Gamma) &\le &
\delta_{\theta}^{\eps_0}(o,o)\,,\quad \mbox{hence}\\
\limsup_{j\to\infty} \delta_{\theta_j}(\Gamma) = 
\liminf_{\eps_0\to 0}\Big(\!\limsup_{j\to\infty} \delta_{\theta_j}(\Gamma)\!\Big) &\le &
\liminf_{\eps_0\to 0}\delta_{\theta}^{\eps_0}(o,o)=\delta_{\theta}(\Gamma)\,.
%\tag*{$\scriptstyle\square$}
\end{eqnarray*}
\qed\\[2mm]

\noindent{\sc Example:}$\quad$
Suppose $\XX$ is a product $\XX=\XX_1\times \XX_2$ of Hadamard manifolds with pinched negative curvature, and 
$\Gamma_1\subset\is(\XX_1)$, $\Gamma_2\subset\is(\XX_2)$  are convex cocompact
groups with critical exponents $\delta_1, \delta_2$. Then by 
Theorem~6.2.5 in \cite{MR1348871} there exists a constant $C>1$ \st for all $n\in\NN $ we have 
\begin{equation}\label{groest}\frac1{C} e^{\delta_i n}\le \#\{\gamma_i\in\Gamma_i : \,n-1<
d(\xo_i,\gamma_i \xo_i)\le n\}\le C e^{\delta_i n}\,,\quad  i=1,2\,.\end{equation}
We are going to examine the action of the product group
$\Gamma=\Gamma_1\times \Gamma_2\subseteq\is(\XX)$ on the product manifold $\XX$. Given $\theta\in (0,\pi/2)$, we estimate for $\eps>0$ sufficiently small the number of orbit points 
\be\Delta
N_\theta^\eps(\xo,\xo;n)&=&\#\{\gamma=(\gamma_1,\gamma_2)\in\Gamma :\,n-1< \sqrt{ d_1(\xo_1,\gamma_1
\xo_1)^2+d_2(\xo_2,\gamma_2\xo_2)^2} \le n\,,\\
&&\qquad |\theta(\xo,\gamma\xo) -\theta|<\eps\}\,\\
&\le&\#\{\gamma=(\gamma_1,\gamma_2)\in\Gamma :\,n-1<\frac{ d_1(\xo_1,\gamma_1 \xo_1)}{\cos\theta(\xo,\gamma\xo)}\le n\,,\\
&&\qquad n-1< \frac{d_2(\xo_2,\gamma_2\xo_2)}{\sin\theta(\xo,\gamma\xo)}\le n\,,\ |\theta(\xo,\gamma\xo) -\theta|<\eps\}\,\\
&\le& C^2 \cdot n\,  e^{\delta_1
  n\cos(\theta-\eps)}\cdot e^{\delta_2n\sin(\theta+\eps)}\,.\ee
As a lower bound, we obtain
\be \Delta N_\theta^\eps(\xo,\xo;n)&\ge&\#\{(\gamma_1,\gamma_2)\in\Gamma : \,n-1< \frac{d_1(\xo_1,\gamma_1\xo_1)}{\cos\theta}\le n\,,\\
&&\qquad n-1<\frac{d_2(\xo_2,\gamma_2\xo_2)}{\sin\theta}\le n\}\ge\frac1{C^2}\cdot e^{\delta_1 n\cos\theta}\cdot e^{\delta_2 n\sin\theta}\,\ee
and therefore conclude $\delta_{\theta}(\Gamma)=\delta_1\cos\theta +\delta_2\sin\theta$. Treating the cases $\theta=0$ and $\theta=\pi/2$ separately one can easily verify that this equation holds for all $\theta\in [0,\pi/2]$.

\section{A generic product for $\Gamma$}

Denote $\RR_{\ge 0}:=\{t\in\RR:t\ge 0\}$. For convenience, we extend the exponent of growth to a map $\Psi_\Gamma:\RR_{\ge 0}^2\to \RR$ as follows: If $x=(x_1,x_2)\in\RR_{\ge 0}^2$ we put $\theta(x):=\arctan (x_2/x_1)$ and set
$$\Psi_\Gamma(x):=||x||\cdot \delta_{\theta(x)}\,.$$
In the remainder of this section we will show that $\Psi_\Gamma$ is a concave function, i.e. for any $x$, $y\in\RR_{\ge 0}^2$ and $t\in [0,1]$ we have
$\Psi_\Gamma(tx+(1-t)y)\ge t\Psi_\Gamma(x)+(1-t)\Psi_\Gamma(y)$. 

Recall that  $\XX$ is a product of locally compact 
Hadamard spaces $\XX_1$, $\XX_2$, $\xo=(\xo_1, \xo_2)$ a fixed base point, and $\Gamma\subset\is(\XX_1)\times\is(\XX_2)$ acts properly discontinuously
and contains a pair of  isometries $g=(g_1,g_2)$ , $h=(h_1,h_2)$  \st $g_i$ and $h_i$ are independent rank one elements in $\Gamma_i$ for $i=1,2$. 
Notice that the distance vector $H:\XX\times \XX\to\RR^2$ defined at the beginning of Section~\ref{prodHadspaces} induces a map 
$\Gamma\to\RR^2$ via the assignment $\gamma\mapsto H(\xo,\gamma\xo)$. By abuse of notation we will call this map also $H$. 

Let  $D$ denote the Dirac measure and $\nu_\Gamma:=\sum_{\gamma\in\Gamma}D_{H(\gamma)}$ the counting measure on $\RR^2$.  In a metric space we denote $B(x,r)$ the ball of radius $r\ge 0$ centered at  $x$.  We will use the following special case of a theorem due to J.-F.~Quint.
\begin{thr}\label{concavity}(\cite{MR1933790}, Theorem 3.2.1)
If there exist $r,s,c>0$ \st for any $x,y\in\RR^2$ the inequality
\begin{equation}\label{concavegrowth}
\nu_\Gamma(B(x+y,s))\ge c\cdot \nu_\Gamma (B(x,r))\cdot \nu_\Gamma (B(y,r))
\end{equation}
 holds, then $\Psi_\Gamma$ is concave.
\end{thr}

In order to prove inequality~(\ref{concavegrowth}) we will construct a generic product for $\Gamma$ as in  \cite{MR1933790}, Proposition 2.3.1. The idea behind is to find a finite set in $\Gamma\times \Gamma$ which maps pairs of orbit points $(\gamma\xo, \varphi^{-1}\xo)$ close to a set $\Ax(g)$ or $\Ax(h)$ as in Definition~\ref{axialisos}.   Unfortunately, unlike in the case of symmetric spaces, we do not dispose of an equivalent of the result of Abels-Margulis-Soifer (\cite[Proposition~2.3.4]{MR1933790}) which plays a crucial role there. Instead, we will exploit the dynamics of a free subgroup of $\langle g, h\rangle\subseteq \Gamma$.

\begin{prp}\label{genprod}
If  $\Gamma\subset \is(\XX_1)\times\is(\XX_2)$ is as above, then there exists a map $\pr:\Gamma\times\Gamma\to \Gamma$ with the following properties:
\begin{enumerate}
\item[(a)] There exists $\kappa\ge 0$ \st for all $\gamma,\varphi\in\Gamma$ we have
$$ \Vert H\big(\pr(\gamma,\varphi)\big)-H(\gamma)-H(\varphi)\Vert\le \kappa\,.$$ 
\item[(b)] For any $r>0$ there exists a finite set $\Lambda\subset \Gamma$ \st for all $\gamma,\varphi,\hat\gamma,\hat\varphi\in \Gamma$ with $\Vert H(\gamma)-H(\hat\gamma)\Vert\le r$, $\Vert H(\varphi)-H(\hat\varphi)\Vert\le r$ we have
$$ \pr(\gamma,\varphi)=\pr(\hat\gamma,\hat\varphi)\quad\Longrightarrow\qquad \hat\gamma\in \gamma \Lambda \ \an \ \hat\varphi\in \Lambda \varphi\,.$$
\end{enumerate}
\end{prp}
\prf \ For $\eta\in \{g^-,g^{+},h^-,h^{+}\}$ and  $i\in\{1,2\}$  let $U_i(\eta)\subset \ganz_i$ be a small neighborhood of $\eta_i$ with $\xo_i\notin U_i(\eta)$ \st all $U_i(\eta)$ are pairwise disjoint. Upon taking smaller neighborhoods, Lemma~\ref{joinrankone} provides a constant  $c>0$  \st for $i\in\{1,2\}$ any pair of points  in distinct neighborhoods can be joined by a rank one geodesic $\sigma_i\subset\XX_i$ with $d(\xo_i,\sigma_i)\le c$.  

In order to construct a map satisfying property~(a) we  let $\gamma=(\gamma_1,\gamma_2)$, $\varphi=(\varphi_1,\varphi_2)\in\Gamma$ arbitrary.
Arguing as in the proof of Proposition~\ref{gettingregularaxials} there exist a finite set $\Lambda\subset\Gamma$ and $\alpha=\alpha(\varphi)$, $\beta=\beta(\gamma)\in\Lambda$ \st
$$\beta\gamma^{-1}\xo\in U_1(h^-)\times U_2(h^-)\ \an\quad \alpha\varphi\xo\in U_1(h^+)\times U_2(h^+)\,.$$
As in the proof of Theorem~\ref{propertiesofPGamma} we set $L:=\max\{d_i(\xo_i,\lambda_i\xo_i): i\in\{1,2\}, \lambda\in\Lambda \}$. For $i=1,2$ we choose a point  $x_i$  on the geodesic joining $\beta_i\gamma_i^{-1}\xo_i$  to $\alpha_i\varphi_i\xo_i$  with $d_i(\xo_i,x_i)\le c$.  Then 
$$d_i(\gamma_i\beta_i^{-1}\alpha_i\varphi_i\xo_i,\xo_i)=d_i(\alpha_i\varphi_i\xo_i,\beta_i\gamma_i^{-1}\xo_i)=d_i(\alpha_i\varphi_i\xo_i,x_i)+d_i(x_i, \beta_i\gamma_i^{-1}\xo_i)\,$$
and we can estimate
\be
d_i(\gamma_i\beta_i^{-1}\alpha_i\varphi_i\xo_i,\xo_i)&\le & d_i(\alpha_i\varphi_i\xo_i,\alpha_i\xo_i)+\overbrace{d_i(\alpha_i\xo_i,\xo_i)}^{\le L} +\overbrace{d_i(\xo_i,x_i)}^{\le c}\\
&&+d_i(x_i,\xo_i)+d_i(\xo_i,\beta_i\xo_i)+d_i(\beta_i\xo_i,\beta_i\gamma_i^{-1}\xo_i)\\
&\le  & d_i(\varphi_i\xo_i,\xo_i)+d_i(\gamma_i\xo_i,\xo_i) +2c +2 L\qquad\mbox{and}\\
d_i(\gamma_i\beta_i^{-1}\alpha_i\varphi_i\xo_i,\xo_i)&\ge & d_i(\varphi_i\xo_i,\xo_i)+d_i(\gamma_i\xo_i,\xo_i) -2c -2 L\,. 
\ee
This gives
$$ \Vert H(\gamma\beta^{-1}\alpha\varphi)-H(\varphi)-H(\gamma)\Vert  \le 2\sqrt2 (c+L)=:\kappa\,,$$
hence the assignment $\pr(\gamma,\varphi):=\gamma\beta(\gamma)^{-1}\alpha(\varphi)\varphi$ satisfies  property (a).

It remains to prove that the map $\pr$ from above also satisfies property (b). Suppose there exists $r>0$ \st  for any finite set $\Lambda_n\subseteq\{\gamma\in\Gamma:\; d(\xo,\gamma\xo)\le n\}$ with $n\in\NN$ there exist $\gamma_n,\varphi_n,\hat\gamma_n,\hat\varphi_n$ with $\Vert H(\gamma_n)-H(\hat\gamma_n)\Vert\le r$, $\Vert H(\varphi_n)-H(\hat\varphi_n)\Vert\le r$ and $g_n:=\pr(\gamma_n,\varphi_n)=\pr(\hat\gamma_n,\hat\varphi_n)$, but $\gamma_n^{-1}\hat\gamma_n\notin \Lambda_n$ or $\hat\varphi_n\varphi_n^{-1}\notin \Lambda_n$. 

Passing to a subsequence if necessary we may assume that all the sequences $(\gamma_n^{-1}\xo)$, $(\hat\gamma_n^{-1}\xo)$, $(\varphi_n\xo)$, $(\hat\varphi_n\xo)\subset\XX$ converge. Notice that even though one of the projections of the sequences to $\XX_1$ or $\XX_2$ may not converge to a boundary point, the arguments from the proof of Proposition~\ref{gettingregularaxials} show that there exist a finite set $\Lambda\subset\Gamma$ and  
$\alpha$, $\hat \alpha$, $\beta$, $\hat\beta\in\Lambda$ \st for all $n\in\NN$
\begin{equation}\label{elementinmenge}
\beta\gamma_n^{-1}\xo,\, \hat\beta\hat\gamma_n^{-1}\xo \in U_1(h^-)\times U_2(h^-) \quad \an\qquad 
\alpha\varphi_n\xo, \, \hat\alpha\hat\varphi_n\xo\in U_1(h^+)\times U_2(h^+) \,.
\end{equation}
For $n\in\NN$ and $i=1,2$ we denote $x_{n,i}$ a point on the geodesic path from $\beta_i\gamma_{n,i}^{-1}\xo_i$ to $\alpha_i\varphi_{n,i}\xo_i$, and $\hat{x}_{n,i}$ a point on the geodesic path from  $\hat\beta_i\hat\gamma_{n,i}^{-1}\xo_i$ to $\hat\alpha_i\hat\varphi_{n,i}\xo_i$ \st $d_i(\xo_i,x_{n,i})\le c$ and $d_i(\xo_i,\hat x_{n,i})\le c$. Furthermore, using  
 $g_n=\gamma_n\beta^{-1}\alpha\varphi_n=\hat\gamma_n\hat\beta^{-1}\hat\alpha \hat\varphi_n$ and  denoting for $i=1,2$ $\sigma_{n,i}$ the geodesic path $\sigma_{\xo_i,g_{n,i}\xo_i}$ there exist $t_i, \hat t_i>0$ \st 
 \be d_i(\gamma_{n,i}\beta^{-1}_i\xo_i, \sigma_{n,i}(t_i))&=& d_i(\gamma_{n,i}\beta_i^{-1}\xo_i, \sigma_{n,i})=d_i(\xo_i, \beta_i\gamma_{n,i}^{-1}\sigma_{n,i})=d_i(\xo_i,x_{n,i})\le c\,,\\
  d_i(\hat \gamma_{n,i}\hat\beta^{-1}_i\xo_i, \sigma_{n,i}(\hat t_i))&=& d_i(\hat\gamma_{n,i}\hat\beta_i^{-1}\xo_i, \sigma_{n,i})=d_i(\xo_i, \hat\beta_i\hat\gamma_{n,i}^{-1}\sigma_{n,i})=d_i(\xo_i,\hat x_{n,i})\le c\,.
 \ee
by~(\ref{elementinmenge}) and Lemma~\ref{joinrankone}. Hence using $L:=\max\{d_i(\xo_i,\lambda_i\xo_i): i\in\{1,2\}, \lambda\in\Lambda \}$
\be d_i(\gamma_{n,i}\xo_i,\sigma_{n,i})&\le & d_i(\gamma_{n,i}\xo_i, \gamma_{n,i}\beta_i^{-1}\xo_i) + d_i(\gamma_{n,i}\beta_i^{-1}\xo_i, \sigma_{n,i}(t_i))\le L+c\,,\\
d_i(\hat \gamma_{n,i}\xo_i,\sigma_{n,i})&\le & d_i(\hat \gamma_{n,i}\xo_i, \hat\gamma_{n,i}\hat\beta_i^{-1}\xo_i) + d_i(\hat \gamma_{n,i}\hat\beta_i^{-1}\xo_i, \sigma_{n,i}(\hat t_i))\le L+c\,.
\ee
For $n\in\NN$ and $i=1,2$ let $y_{n,i}, \hat y_{n,i}\in\XX_i$ be the points on the geodesic path $\sigma_{n,i}$ \st $d_i(\xo_i, y_{n,i})=d_i(\xo_i,\gamma_{n,i}\xo_i)$ and $d_i(\xo_i, \hat y_{n,i})=d_i(\xo_i,\hat\gamma_{n,i}\xo_i)$. Since $\Vert H( \gamma_n)-H(\hat\gamma_n)\Vert\le r$ we have $d_i(y_{n,i},\hat y_{n,i})\le r$, and, by elementary geometric estimates, 
$$d_i(\gamma_{n,i}\xo_i,y_{n,i})\le 2(L+c)\quad\an\qquad d_i(\hat\gamma_{n,i}\xo_i,\hat y_{n,i})\le 2(L+c)\,.$$ 
We summarize 
\be d_i(\xo_i,\gamma_{n,i}^{-1}\hat \gamma_{n,i}\xo_i)&= & d_i(\gamma_{n,i}\xo_i,\hat\gamma_{n,i}\xo_i)\le d_i(\gamma_{n,i}\xo_i,y_{n,i})+d(y_{n,i},\hat y_{n,i})+d_i(\hat y_{n,i},\hat\gamma_{n,i}\xo_i)\\
&\le &  2(L+c) + r + 2(L+c)\,,\ \mbox{i.e.}\\
d(\xo,\gamma_n^{-1}\hat \gamma_n\xo)&\le & \sqrt{2}(4L+4c+r)=:R\,.
\ee
In particular, for $n>R$ we have $\gamma_n^{-1}\hat\gamma_n\in \Lambda_n$, and, in order to obtain the desired contradiction, it remains to prove that $\hat\varphi_n\varphi_n^{-1}\in \Lambda_n$ for $n$ sufficiently large. 

Notice that $\hat\varphi_n=\hat\alpha^{-1}\hat\beta\hat\gamma_n^{-1}g_n=\hat\alpha^{-1}\hat\beta\hat\gamma_n^{-1}\gamma_n\beta^{-1}\alpha\varphi_n$, hence
\begin{eqnarray*}
d(\xo,\hat\varphi_n\varphi_n^{-1}\xo)&=& d(\xo, \hat\alpha^{-1}\hat\beta\hat\gamma_n^{-1}\gamma_n\beta^{-1}\alpha\xo)\le \overbrace{d(\xo,\hat\alpha^{-1}\xo)}^{\le\sqrt{2}L}+\overbrace{d(\hat\alpha^{-1}\xo,\hat\alpha^{-1}\hat\beta\xo)}^{\le\sqrt{2}L}+\\
&&d(\hat\alpha^{-1}\hat\beta\xo,\hat\alpha^{-1}\hat\beta\hat\gamma_n^{-1}\gamma_n\xo)+ d(\hat\gamma_n^{-1}\gamma_n\xo,\hat\gamma_n^{-1}\gamma_n\beta^{-1}\xo)+d(\beta^{-1}\xo,\beta^{-1}\alpha\xo)\\
&\le&d(\xo,\gamma_n^{-1}\hat\gamma_n\xo) + 4\sqrt{2}L\le R+ 4\sqrt{2}L\,. \end{eqnarray*}
This finishes the proof. \qed\\[3mm]
The following lemma now shows that equation~(\ref{concavegrowth}) holds. 
\begin{lem}
There exist $r,s,c>0$ \st for any $x,y\in\RR^2$ we have
$$\nu_\Gamma(B(x+y,s))\ge c\cdot \nu_\Gamma (B(x,r))\cdot \nu_\Gamma (B(y,r))\,.$$
\end{lem}
\prf\  Notice that $\nu_\Gamma (B(x,r))=\#\{\gamma\in\Gamma:\Vert H(\gamma)-x\Vert <s\}$. Fix $r>0$, put $s=\kappa+2r$ with $\kappa\ge 0$ from  Proposition~\ref{genprod} (a) and denote $C>0$ the inverse of the cardinality of the set $\Lambda\times\Lambda$ from  Proposition~\ref{genprod} (b). Put 
$$P(\Gamma):=\{(\gamma,\varphi)\in\Gamma\times\Gamma :\Vert H(\gamma)-x\Vert <r\,,\ \Vert H(\varphi)-y\Vert <r\}\,.$$
We will show that for all $x,y\in\RR^2$ 
$$\#\{\gamma\in\Gamma : \Vert H(\gamma)-x-y\Vert<s\}\ge C\cdot \# P(\Gamma)\,.$$
Let $(\gamma,\varphi)\in P(\Gamma)$. Then $\alpha:=\pr(\gamma,\varphi)\in\Gamma$ satisfies
\be \Vert H(\alpha)-x-y\Vert &\le &\Vert H(\alpha)-H(\gamma)-H(\varphi)\Vert +\Vert H(\gamma)-x\Vert+\Vert H(\varphi)-y\Vert\\
&\le& \kappa+r+r=s\,.\ee
Moreover,  Proposition~\ref{genprod} (b) implies that the number of  different elements in $P(\Gamma)$ which can yield the same element in $\{\gamma\in\Gamma : \Vert H(\gamma)-x-y\Vert<s\}$ is bounded by $\#(\Lambda\times\Lambda)$.\qed\\[3mm]
As a corollary of Theorem~\ref{concavity} and Proposition~\ref{genprod} we obtain
\begin{thr}
The function $\Psi_\Gamma$ is concave.
\end{thr}
Together with Proposition~\ref{upsemcon} this gives Theorem~E of the introduction.

\bibliography{Bibliographie} 

\end{document}